\documentclass[11pt, DIV10,a4paper]{article}

\usepackage{natbib}
\newcommand {\ctn}{\citet} 
\usepackage{graphicx,subfigure,amsmath,latexsym,amssymb}
\usepackage{float,epsfig,multirow,rotating,times}
\usepackage{upgreek,wrapfig}
\usepackage{comment}

\newcommand{\btheta}{\boldsymbol{\theta}}

\newcommand{\bbeta}{\boldsymbol{\beta}}

\newcommand{\bxi}{\boldsymbol{\xi}}
\newcommand{\bGamma}{\boldsymbol{\Gamma}}
\newcommand{\bTheta}{\boldsymbol{\Theta}}
\newcommand{\bgamma}{\boldsymbol{\gamma}}

\newcommand{\bvarsigma}{\boldsymbol{\varsigma}}

\newcommand{\bC}{\boldsymbol{C}}

\newcommand{\bI}{\boldsymbol{I}}

\newcommand{\bZ}{\boldsymbol{Z}}
\newcommand{\bz}{\boldsymbol{z}}

\newtheorem{theorem}{Theorem}

\newtheorem{corollary}[theorem]{Corollary}

\newtheorem{lemma}[theorem]{Lemma}

\newtheorem{remark}[theorem]{Remark}

\newenvironment{proof}[1][Proof]{\textbf{#1.} }{\ \rule{0.5em}{0.5em}}

\newcommand{\boldcal}[1]{\mbox{\boldmath{$\mathcal #1 $}}}

\numberwithin{equation}{section}
\numberwithin{algo}{section}
\numberwithin{table}{section}
\numberwithin{figure}{section}

\usepackage{setspace}
\usepackage[mathscr]{euscript}
\usepackage[margin=1in]{geometry}
\begin{document}

\normalsize

\title{\vspace{-0.8in}
Asymptotic Theory of Bayes Factor in Stochastic Differential Equations: Part I}
\author{Trisha Maitra and Sourabh Bhattacharya\thanks{
Trisha Maitra is a PhD student and Sourabh Bhattacharya 
is an Associate Professor in
 Interdisciplinary Statistical Research Unit, Indian Statistical
Institute, 203, B. T. Road, Kolkata 700108.
Corresponding e-mail: sourabh@isical.ac.in.}}
\date{\vspace{-0.5in}}
\maketitle%

\begin{abstract}
Research on asymptotic model selection in the context of stochastic differential equations ($SDE$'s) is almost non-existent
in the literature. In particular, when a collection of $SDE$'s is considered, 
the problem of asymptotic model selection has not been hitherto investigated. 
Indeed, even though the diffusion coefficients may be considered known, questions on appropriate choice of the drift functions
constitute a non-trivial model selection problem.

In this article, we develop the asymptotic theory for comparisons between collections of $SDE$'s with respect
to the choice of drift functions using Bayes factors when 
the number of equations (individuals)
in the collection of $SDE$'s tend to infinity while the time domains remain bounded for each equation.
Our asymptotic theory covers situations when the observed processes associated with the
$SDE$'s are independently and identically distributed ($iid$), as well as when they are independently 
but not identically distributed (non-$iid$). In particular, we allow incorporation of available time-dependent 
covariate information into each $SDE$ through a multiplicative factor of the drift function; 
we also permit different initial values and domains of observations for the $SDE$'s.

Our model selection problem thus encompasses selection of a set of appropriate time-dependent covariates
from a set of available time-dependent covariates, besides selection
of the part of the drift function free of covariates. 

For both $iid$ and non-$iid$ set-ups we establish almost sure exponential convergence
of the Bayes factor. 

Furthermore, we demonstrate with simulation studies that even in non-asymptotic scenarios
Bayes factor successfully captures the right set of covariates.
\\[2mm]
 
\end{abstract}

\section{Introduction}
\label{sec:intro}

In statistical applications where ``within" subject variability is caused by some random component
varying continuously in time, stochastic differential equations ($SDE$'s)
have important roles to play for modeling the temporal component of each individual.
The inferential abilities of the $SDE$'s can be enhanced by incorporating
covariate information available for the subjects. In these time-dependent situations
it is only natural that the available covariates are also continuously varying with time. 
Examples of statistical applications of $SDE$-based models with time-dependent covariates are
\ctn{Zita11}, \ctn{Overgaard05}, \ctn{Leander15}, the first one also considering the hierarchical Bayesian paradigm.

Unfortunately, asymptotic inference in systems of $SDE$ based models consisting of time-varying covariates seem to be
rare in the statistical literature, in spite of their importance. 
So far random effects $SDE$ models have been considered for asymptotic inference, without covariates.
We refer to \ctn{Maud12} for a brief review, who also undertake theoretical and classical asymptotic 
investigation of a class of
random effects models based on $SDE$'s. Specifically, they model the $i$-th individual by
\begin{equation}
d X_i(t)=b(X_i(t),\phi_i)dt+\sigma(X_i(t))dW_i(t),
\label{eq:sde_basic1}
\end{equation}
where, for $i=1,\ldots,n$, $X_i(0)=x^i$ is the initial value of the stochastic process $X_i(t)$, which 
is assumed to be continuously observed on the time interval $[0,T_i]$; $T_i>0$ assumed to be known.
The function $b(x,\varphi)$, which is the drift function, is a known, real-valued function on $\mathbb R\times\mathbb R^d$ 
($\mathbb R$ is the real line and $d$ is the dimension), and the function $\sigma:\mathbb R\mapsto\mathbb R$ is the known 
diffusion coefficient.
The $SDE$'s given by (\ref{eq:sde_basic1}) are driven by independent standard Wiener processes $\{W_i(\cdot);~i=1,\ldots,n\}$, 
and $\{\phi_i;~i=1,\ldots,n\}$, which are to be interpreted as the random effect parameters associated
with the $n$ individuals, which are assumed by \ctn{Maud12} to be independent of the Brownian motions and
independently and identically distributed ($iid$) random variables with some common distribution.
For the sake of convenience 
\ctn{Maud12} (see also \ctn{Maitra14a} and \ctn{Maitra14b}) assume $b(x,\phi_i)=\phi_ib(x)$. 
Thus, the random effect is a multiplicative factor of the drift function. 
In this article, we generalize the multiplicative factor to include time-dependent covariates.

In the case of $SDE$-based models,   
proper specification of the drift function and the associated prior distributions demand
serious attention, and this falls within the purview of model selection. 
Moreover, when (time-varying) covariate information are available, there arises the problem
of variable selection, that is, the most appropriate subset from the set of many available covariates needs to be chosen. 
As is well-known (see, for example, \ctn{Kass95}), the Bayes factor (\ctn{Jeffreys61}) is a strong candidate for dealing with
complex model selection problems. Hence, it is natural to consider this criterion for model selection
in $SDE$ set-ups. However, dealing with Bayes factors directly in $SDE$ set-ups is usually infeasible
due to unavailability of closed form expressions, and hence various numerical approximations based
on Markov chain Monte Carlo, as well as related criteria such
as Akaike Information Criterion (\ctn{Akaike73}) and Bayes Information Criterion (\ctn{Schwarz78}), are generally employed
(see, for example, \ctn{Fuchs13}, \ctn{Iacus08}).  
But quite importantly, although Bayes factor and its variations find use in general $SDE$ models, 
in our knowledge covariate selection 
in $SDE$ set-ups has not been addressed so far. 

Moreover, asymptotic theory of Bayes factors
in $SDE$ contexts, with or without covariates, is still lacking 
(but see \ctn{Siva02} who asymptotically compare three specific diffusion models in
single equation set-ups using intrinsic and fractional Bayes factors).  
In this paper, our goal is to develop an asymptotic theory of Bayes factors
for comparing different sets of $SDE$ models. Our asymptotic theory simultaneously
involves time-dependent covariate selection associated with a multiplicative
part of the drift function, in addition to selection of the part of the drift
function free of covariates. The asymptotic framework of this paper assumes
that the number of individuals tends to infinity, while their domains of observations
remain bounded.
%

It is important to clarify that the diffusion coefficient is not associated
with the question of model selection. Indeed, it is already known from \ctn{Robert01} 
that when the associated continuous process is 
completely observed, the diffusion coefficient of the relevant $SDE$ can be calculated directly. Moreover, 
two diffusion processes with different diffusion coefficients are orthogonal. 
Consequently, we assume throughout that the diffusion coefficient of the $SDE$'s is known.


We first develop the model selection theory using Bayes factor in general $SDE$ 
based $iid$ set-up; note that the $iid$ set-up
ensues when there is no covariate associated with the model and when the initial values and the domains of observations
are the same for every individual. The model selection problem in $iid$ cases is essentially associated 
with the choice of the drift
functions with no involvement of covariate selection.
We then extend our theory to the 
non-$iid$ set-up, consisting of time-varying covariates and different  
initial values and domains of observations. Here model selection involves not only selection of 
the part of the drift functions free of the covariates, but also the subset of important covariates
from a set of available covariates.

Specifically, we prove almost sure exponential convergence of the relevant Bayes
factors in our set-ups. 
Assuming the $iid$ set-up we develop our asymptotic
theory based on a general result already existing in the literature.
However, for the non-$iid$ situation we first develop a general theorem which may perhaps
be of independent interest, 
and prove almost sure exponential convergence of the Bayes factor in our non-$iid$ $SDE$ set-up as a special case
of our theorem. 

It is important to note that (which we also clarify subsequently in Section \ref{subsec:no_random_effects}), 
that in the asymptotic framework of this paper, where
the domains of observations remain bounded for the individuals, incorporation of random effects does not make sense
from the asymptotic perspective. For this reason we include random effects in our paper \ctn{Maitra15b}, where
we assume that even the domains of observations are allowed to increase indefinitely.

The rest of our article is structured as follows. In Section \ref{sec:bf_sde} we formalize the problem of model
selection in our aforementioned asymptotic framework. 
We then present the necessary assumptions and results in Section \ref{sec:case1}. In Section \ref{sec:asymp_BF_iid} 
we investigate convergence of the Bayes factor when the $SDE$ models being compared form an 
$iid$ system of equations. In Section \ref{sec:non_iid} we develop a general asymptotic theory of 
Bayes factors in the non-$iid$ situation, and then in Section \ref{sec:sde_non_iid} we investigate 
exponential convergence of the Bayes factor when the system of $SDE$'s are non-$iid$. 
In Section \ref{sec:simulated_data} we demonstrate with simulation studies that Bayes factor
yields the correct covariate combination for our $SDE$ models even in non-asymptotic cases.
We provide a brief summary of this article and make concluding remarks in Section \ref{sec:conclusion}. 

The proofs of our lemmas and theorems are provided in the supplementary document 
whose sections will be referred to in this article by the prefix ``S-".

\section{Formalization of the model selection problem in the $SDE$ set-up}
\label{sec:bf_sde}

Our assumptions (H2$^\prime$)  in Section \ref{sec:case1} 
ensure that our considered systems are well defined and we are able to compute the exact likelihood. 
We consider the filtration ($\mathcal F_t^W,t\geq 0$),
where $\mathcal F_t^W=\sigma(W_i(s),s\leq t)$. Each process $W_i$ is a $(\mathcal F_t^W , t\geq 0)$-adapted Brownian
motion.

In connection with model selection we must analyze the same data set with respect to two different models. 
So, although the distribution of the underlying stochastic process under the two models are different, 
to avoid notational complexity we denote the process by $X_i(t)$ under both the models, keeping in mind that
the distinction becomes clear from the context and also by the model-specific parameters.

\subsection{The structure of the $SDE$ models to be compared}
\label{subsec:basic_models}
Now, let us consider the following two systems of $SDE$ models for $i=1,2,\ldots,n$:
\begin{equation}
d X_i(t)=\phi_{i,\bxi_0}(t)b_{\bbeta_0}(X_i(t))dt+\sigma(X_i(t))dW_i(t)
\label{eq:sde1}
\end{equation}
and
\begin{equation}
d X_i(t)=\phi_{i,\bxi_1}(t)b_{\bbeta_1}(X_i(t))dt+\sigma(X_i(t))dW_i(t)
\label{eq:sde2}
\end{equation}
where, $X_i(0)=x^i$ is the initial value of the stochastic process $X_i(t)$, which is assumed to be continuously 
observed on the time interval $[0,T_i]$; $T_i>0$ for all $i$ and assumed to be known.
We assume that (\ref{eq:sde1}) represents the true model and (\ref{eq:sde2}) is any other model.
In the above equations, for $j=0,1$, $\bxi_j$ and $\bbeta_j$ denote the sets of parameters
associated with the true model and the other model.

\subsection{Incorporation of time-varying covariates}
\label{subsec:covariates}
We model $\phi_{i,\bxi_j}(t)$ for $j=0,1$, as  
\begin{equation}
\phi_{i,\bxi_j}(t)=\phi_{i,\bxi_j}(\bz_i(t))
=\xi_{0j}+\xi_{1j}g_1(z_{i1}(t))+\xi_{2j}g_2(z_{i2}(t))+\cdots+\xi_{pj}g_p(z_{ip}(t)),
\label{eq:phi_model}
\end{equation}
where $\bxi_j=(\xi_{0j},\xi_{1j},\ldots,\xi_{pj})$ is a set of real constants for $j=0,1$, 
and $\bz_i(t)=(z_{i1}(t),z_{i2}(t),\ldots,z_{ip}(t))$ is the set of available covariate information 
corresponding to the $i$-th individual, depending upon time $t$. We assume $\bz_i(t)$ is continuous in $t$,
$z_{il}(t)\in \bZ_l$  where $\bZ_l$ is  compact and 
$g_l:\bZ_l\rightarrow \mathbb R$ is  continuous, for $l=1,\ldots,p$.
We let $\boldcal Z=\bZ_1\times\cdots\times\bZ_p$, and 
$\mathfrak Z=\left\{\bz(t)\in\boldcal Z:t\in[0,\infty)~\mbox{such that}~\bz(t)~\mbox{is continuous in}~t\right\}$.
Hence, $\bz_i\in\mathfrak Z$ for all $i$.
The functions $b_{\bbeta_j}$ are multiplicative parts of the drift functions free of the covariates.

\subsection{Model selection with respect to the drift function and the covariates}
\label{subsec:model_selection_drift_covariates}

We accommodate the possibility that the dimensions of $\bbeta_0,\bbeta_1$, associated with the drift functions, may be different.
In reality, $b_{\bbeta_0}$ may be piecewise linear or convex combinations of linear functions, where
the number of linear functions involved (and hence, the number of associated intercept and slope parameters)
may be unknown. That is, not only the values of the components of the parameter $\bbeta_0$, but also the 
number of the components of $\bbeta_0$ may be unknown in reality. 
In general, $b_{\bbeta_0}$ may be any function, linear or non-linear, satisfying some desirable conditions.
Linearity assumptions may be convenient, but need not necessarily be unquestionable.
In other words, modeling $b_{\bbeta_0}$ in the $SDE$ context is a challenging exercise,
and hence the issue of model selection in this context must play an important role in the $SDE$ set-up.

We also accommodate the possibility that $\bxi_0$ and $\bxi_1$, associated with $\phi_{i,\bxi_0}$ and $\phi_{i,\bxi_1}$,
may be coefficients associated with different subsets of the available set of $p$ covariates. 
This has important implication from the viewpoint of variable selection.
Indeed, in a set of $p$ time-dependent covariates, all the covariates are unlikely to be significant, particularly
if $p$ is large.
Thus, some (perhaps, many) of the coefficients $\xi_{l0}$ associated with the true model must be zero. 
This means that only a specific subset of the $p$ covariates is associated with the true model. If a different
set of covariates, associated with $\bxi_1$, is selected for actually modeling the data, then the Bayes factor
is expected to favour the true set of covariates associated with $\bxi_0$. 

If two different models are compared by the Bayes factor, none of which may be the true model, then the
Bayes factor is expected to favour that model which is closest to the true model in terms of the Kullback-Leibler
divergence.

\subsection{Form of the Bayes factor}
\label{subsec:bf_form}
For $j=0,1$, letting $\btheta_j=(\bbeta_j,\bxi_j)$, we
first define the following quantities:
\begin{equation}
U_{i,\btheta_j} =\int_0^{T_i}\frac{\phi_{i,\bxi_j}(s)b_{\bbeta_j}(X_i(s))}{\sigma^2(X_i(s))}dX_i(s),
\quad\quad V_{i,\btheta_j} =\int_0^{T_i}\frac{\phi^2_{i,\bxi_j}(s)b_{\bbeta_j}^2(X_i(s))}{\sigma^2(X_i(s))}ds
\label{eq:u_v}
\end{equation}
for $j=0,1$ and $i=1,\ldots,n$. 

Let $\bC_{T_i}$ denote the space of real continuous functions $(x(t), t \in [0,T_i ])$ defined on $[0,T_i ]$, 
endowed with the $\sigma$-field $\mathcal C_{T_i}$ associated with the topology of uniform convergence 
on $[0,T_i ]$. We consider the distribution $P^{x_i,T_i,\bz_i}_j$ on $(C_{T_i} ,\mathcal C_{T_i})$ 
of $(X_i (t), t\in [0,T_i ])$ given by (\ref{eq:sde1}) and (\ref{eq:sde2}) for $j=0,1$. 
We choose the dominating measure $P_i $ as the distribution of (\ref{eq:sde1}) and (\ref{eq:sde2}) with null
drift. So, for $j=0,1$, 
\begin{align}
\frac{dP^{x_i,T_i,\bz_i}_j}{dP_i}=f_{i,\btheta_j}(X_i)
&=\exp\left(U_{i,\btheta_j}-\frac{V_{i,\btheta_j}}{2}\right), 
\label{eq:densities}
\end{align}
where $f_{i,\btheta_0}(X_i)$ denotes the true density and $f_{i,\btheta_1}(X_i)$ stands for the other density
associated with the modeled $SDE$.

Let $\bTheta=\mathfrak B\times\bGamma$ be the parameter space on which a prior probability 
measure of $\btheta_1$, which we denote by $\pi(\btheta_1)$, is proposed. 
In the set-up where $n\rightarrow\infty$ and $T_i$ are given, 
we are interested in asymptotic properties of the Bayes factor, given by,
$I_0\equiv 1$ and for $n\geq 1$,
\begin{align}
I_n&=\int_{\bTheta} R_n(\btheta_1)\pi(d\btheta_1),
\label{eq:I_n}
\end{align}
as $n\rightarrow\infty$, where
\[
R_n(\btheta_1)=\prod_{i=1}^n\frac{f_{i,\btheta_1}(X_i)}{f_{i,\btheta_0}(X_i)}.
\]

\subsection{The $iid$ and the non-$iid$ set-ups}
\label{subsec:iid_non_iid}

Note that, for $iid$ set-up $\btheta_j=(\bbeta_j,\xi_{0j})$, along with $x^i=x$ and $T_i=T$ for all $i$.
Since, for the $iid$ set-up $\bxi_j=\xi_{0j}$, so, in this case 
$\bGamma=\mathbb R$. 
Thus, here the problem of model selection reduces to comparing $\xi_{00}b_{\bbeta_0}$ with $\xi_{01}b_{\bbeta_1}$
using Bayes factor.

In the non-$iid$ set-up we relax the assumptions $\xi_{1j}=\xi_{2j}=\cdots=\xi_{pj}=0$ and
$x^i=x$, $T_i=T$ for each $i$. Hence, in this case, the model selection problem involves
variable selection as well as comparison between different drift functions.

\subsection{No random effects when $T_i$ are given}
\label{subsec:no_random_effects}
It is important to perceive that when the $T_i$ are fixed constants, it is not possible to
allow random effects into the model and still achieve consistency of the Bayes factor. This is because
in that case the $SDE$ set-up would simply reduce to $n$ independent models, each with independent sets
of parameters, leaving no scope for asymptotics since $T_i$ are held constants.
In \ctn{Maitra15b} we consider random effects when $T_i\rightarrow\infty$ along with $n\rightarrow\infty$.

\subsection{A key relation between $U_{i,\btheta_j}$ and $V_{i,\btheta_j}$ in the context
of model selection using Bayes factors}
\label{subsec:key_relation_U_V}
An useful relation between $U_{i,\btheta_j}$ and $V_{i,\btheta_j}$ which we will often make use of
in this paper is as follows.
\begin{align}
U_{i,\btheta_j}&=\int_0^{T_i}\frac{\phi_{i,\bxi_j}(s)b_{\bbeta_j}\left(X_i(s)\right)}
{\sigma^2\left(X_i(s)\right)}dX_i(s)\notag\\
&=\int_0^{T_i}\frac{\phi_{i,\bxi_j}(s)b_{\bbeta_j}\left(X_i(s)\right)}{\sigma^2\left(X_i(s)\right)}
\left[\phi_{i,\bxi_0}(s)b_{\bbeta_0}\left(X_i(s)\right)ds+\sigma\left(X_i(s)\right)dW_i(s)\right]\notag\\
&=\int_0^{T_i}\frac{\phi_{i,\bxi_j}(s)\phi_{i,\bxi_0}(s)b_{\bbeta_j}\left(X_i(s)\right)b_{\bbeta_0}\left(X_i(s)\right)}
{\sigma^2\left(X_i(s)\right)}ds
+\int_0^{T_i}\frac{\phi_{i,\bxi_j}(s)b_{\bbeta_j}\left(X_i(s)\right)}
{\sigma\left(X_i(s)\right)}dW_i(s)\notag\\
&=V_{i,\btheta_0,\btheta_j}
+\int_0^{T_i}\frac{\phi_{i,\bxi_j}(s)b_{\bbeta_j}\left(X_i(s)\right)}
{\sigma\left(X_i(s)\right)}dW_i(s),
\label{eq:u_v_relation}
\end{align}
with
\begin{equation}
V_{i,\btheta_0,\btheta_j}=
\int_0^{T_i}\frac{\phi_{i,\bxi_j}(s)\phi_{i,\bxi_0}(s)b_{\bbeta_j}\left(X_i(s)\right)b_{\bbeta_0}
\left(X_i(s)\right)}{\sigma^2\left(X_i(s)\right)}ds.
\label{eq:V_0_j}
\end{equation}
Note that $V_{i,\btheta_0}=V_{i,\btheta_0,\btheta_0}$ and 
$V_{i,\btheta_1}=V_{i,\btheta_1,\btheta_1}$.
Also note that, for $j=0,1$, for each $i$,
\begin{equation}
E_{\btheta_0}\left[\int_0^{T_i}\frac{\phi_{i,\bxi_j}(s)b_{\bbeta_j}\left(X_i(s)\right)}
{\sigma\left(X_i(s)\right)}dW_i(s)\right]=0,
\label{eq:zero_mean}
\end{equation}
so that $E_{\btheta_0}\left(U_{i,\btheta_j}\right)=E_{\btheta_0}\left(V_{i,\btheta_0,\btheta_j}\right)$.

\section{Requisite assumptions and results for the asymptotic theory of Bayes factor 
when $n\rightarrow\infty$ but $T_i$ are constants for every $i$}
\label{sec:case1}



We assume the following conditions:
\begin{itemize}
\item[(H1$^\prime$)] The parameter space $\bTheta=\mathfrak B\times\bGamma$ such that $\mathfrak B$ 
and $\bGamma$ are compact.
\end{itemize}

\begin{itemize}
\item[(H2$^\prime$)] 
For $j=0,1$, 
$b_{\bbeta_j}(\cdot)$ and  $\sigma(\cdot)$ are $C^1$ on $\mathbb R$ and satisfy
$b^2_{\bbeta_j}(x)\leq K_1(1+x^2+\|\bbeta_j\|^2)$  
and $\sigma^2(x)\leq K_2(1+x^2)$
for all $x\in\mathbb R$, for some $K_1, K_2>0$. Now, due to (H1$^\prime$) 
the latter boils down to assuming
$b^2_{\bbeta_j}(x)\leq K(1+x^2)$  
and $\sigma^2(x)\leq K(1+x^2)$
for all $x\in\mathbb R$, for some $K>0$. 
\end{itemize}
Because of (H2$^\prime$) it follows from Theorem 4.4 of \ctn{Mao11}, page 61, that  
for all $T_i>0$, and any $k\geq 2$, 
\begin{equation}
E\left(\underset{s\in [0,T_i]}{\sup}~|X_i(s)|^k\right)\leq\left(1+3^{k-1}E|X_i(0)|^k\right)\exp\left(\tilde\vartheta T_i\right),
\label{eq:moment1}
\end{equation}
where
$$\tilde\vartheta=\frac{1}{6}\left(18K\right)^{\frac{k}{2}}T^{\frac{k-2}{2}}_i\left[T^{\frac{k}{2}}_i
+\left(\frac{k^3}{2(k-1)}\right)^{\frac{k}{2}}\right].$$
We further assume:
\begin{itemize}
\item[(H3$^\prime$)] For every $x$, $b_{\bbeta_j}(x)$ is continuous in $\bbeta_j$, for $j=0,1$.

\item[(H4$^\prime$)] For $j=0,1$, 
\begin{equation}
\frac{b^2_{\bbeta_j}(x)}{\sigma^2(x)}
\leq K_{\bbeta_j} \left(1+x^2+\|\bbeta_j\|^2\right),
\label{eq:H4_prime_1}
\end{equation}
where $K_{\bbeta_j}$ is continuous in $\bbeta_j$.

\item[(H5$^\prime$)]
(i) Let 
$\boldcal Z= \bZ_1\times\bZ_2\times\cdots\times\bZ_p$ be the space of covariates where $\bZ_l$ is  compact for 
$l=1,\ldots,p$ and $\bz_i(t)=(z_{i1}(t),z_{i2}(t),\ldots,z_{ip}(t))\in\boldcal Z$ for every $i=1,\ldots,n$
and $t\in[0,T_i]$. Moreover, $\bz_i(t)$ are continuous in $t$, so that $\bz_i\in\mathfrak Z$ for every $i$.

(ii) For $j=0,1$, the vector of covariates $\bz_i(t)$ is related to the $i$-th $SDE$
of the $j$-th model via
$$\phi_{i,\bxi_j}(t)=\phi_{\bxi_j}(\bz_i(t))=\xi_{0j}+\sum_{l=1}^p\xi_{lj}g_l(\bz_i(t)),$$
where, for $l=1,\ldots,p$, $g_l:\bZ_l\rightarrow \mathbb R$ is 
continuous. 
Notationally, for a given $\bz(t)$, we denote $\phi_{\bxi_j}(t)=\phi_{\bxi_j}(\bz(t))=\xi_{0j}+\sum_{l=1}^p\xi_{lj}g_l(\bz(t))$.

(iii) For $l=1,\ldots,p$, and for $t\in [0,T_i]$,
\begin{equation}
\frac{1}{n}\sum_{i=1}^n g_l(z_{il}(t))\rightarrow c_{l}(t);
\label{eq:H5_prime_1}
\end{equation}
and, for $l,m=1,\ldots,p$; $t\in [0,T_i]$,
\begin{equation}
\frac{1}{n}\sum_{i=1}^n g_l(z_{il}(t))g_m(z_{im}(t))\rightarrow c_l(t)c_m(t),
\label{eq:H5_prime_2}
\end{equation}
as $n\rightarrow\infty$, where $c_l(t)$ are real constants.
\end{itemize}

Note that, given $l$ and $t$, had $z_{il}(t)$ been random and $iid$ with respect to $i$, then
(\ref{eq:H5_prime_1}) would hold almost surely by the strong law of large numbers. Additionally,
if $z_{il}(t)$ and $z_{im}(t)$ were independent, then (\ref{eq:H5_prime_2}) would hold almost surely
as well. Hence, in this paper, one may assume that for $i=1,\ldots,n$, and 
$l=1,\ldots,p$, the covariates $z_{il}$ are observed realizations
of stochastic processes that are $iid$ for $i=1,\ldots,n$, for all $l=1,\ldots,p$, and that
for $l\neq m$, the processes generating $z_{il}$ and $z_{im}$ are independent. 
Thus, in essence, we assume here that for $l\neq m$, $g_l(z_{il}(t))$ and $g_m(z_{im}(t))$ are 
uncorrelated. 

We then have the following lemma, which will be useful for proving our main results. 

\begin{lemma}
\label{lemma:moment_existence}
Assume (H1$^\prime$) -- (H4$^\prime$). Then for all $\btheta_1\in\mathfrak B\times\bGamma$, for $k\geq 1$,  
\begin{align}
E_{\btheta_0}\left[U_{i,\btheta_j}\right]^{k}&<\infty; j=0,1,\label{eq:finite_moment_u}\\
E_{\btheta_0}\left[V_{i,\btheta_1}\right]^{k}&<\infty,\label{eq:finite_moment_v1}\\
E_{\btheta_0}\left[V_{i,\btheta_0,\btheta_j}\right]^{k}&<\infty; j=0,1.\label{eq:finite_moment_v}
\end{align}
Moreover, for $j=1$, the above expectations are continuous in $\btheta_1$.
\end{lemma}

\section{Convergence of Bayes factor in the $SDE$ based $iid$ set-up}
\label{sec:asymp_BF_iid}

We first consider the $iid$ set-up; in other words, we assume that $x^i=x$, $T_i=T$ 
for $i=1,\ldots,n$, and $j=0,1$. In this case $\btheta_j=(\bbeta_j,\xi_{0j})$ for $j=0,1$. 
We shall relax these assumptions subsequently when we take up the non-$iid$ (that is,
independent, but non-identical) case.



\subsection{A general result on consistency of Bayes factor in the $iid$ set-up}
\label{subsec:result_BF_consistency}

To investigate consistency of the Bayes factor, we resort to a general result 
in the $iid$ set-up developed by \ctn{Walker04a} (see also \ctn{Walker04b}). 
To state the
result we first define some relevant notation
which apply to both parametric and nonparametric problems.
For any $x$ in the appropriate domain, let 
\[ \hat f_n(x)=\int f(x)\pi_n(df)\]
be the posterior predictive density, where $\pi_n$ stands for the posterior of $f$, 
given by
\[\pi_n(A) = \frac{\int_A\prod_{i=1}^nf(X_i)\pi(df)}{\int\prod_{i=1}^nf(X_i)\pi(df)}\]
and let
\[ \hat f_{nA}(x)=\int f(x)\pi_{nA}(df)\]
be the posterior predictive density restricted to the set $A$, that is, for the prior probability $\pi(A)>0$,
\[ \pi_{nA}(df)=\frac{\bI_A(f)\pi_n(df)}{\int_A\pi_n(df)},\]
where $\bI_A$ denotes the indicator function of the set $A$. 

Clearly, the above set-up is in accordance with the $iid$ situation.
The following theorem of \ctn{Walker04a} is appropriate for our $iid$ set-up.

\begin{theorem}[\ctn{Walker04a}]
\label{theorem:walker}
Let $f_0$ be the density of the true data-generating distribution and $f$ be the density of the modeled distribution.
Also let $\mathcal K(f_0,f)=\int f_0(x)\log\left(\frac{f_0(x)}{f(x)}\right)dP_0$ denote the Kullback-Leibler divergence between $f_0$ and $f$, where $P_0$ is the appropriate
dominating measure associated with $f_0$.
Assume that 
\begin{equation}
\pi\left(f:\mathcal K(f_0,f)<c_1\right)>0, 
\label{eq:kl_property}
\end{equation}
only for, and for all $c_1>\delta$, for some $\delta\geq 0$,
and that for all $\epsilon>0$, 
\begin{equation}
\underset{n}{\lim\inf}~\mathcal K\left(f_0,\hat f_{nA(\epsilon)}\right)\geq\epsilon,
\label{eq:Q_star}
\end{equation}
when $A(\epsilon)=\left\{f:\mathcal K\left(f_0,f\right)>\epsilon\right\}$. 
Property (\ref{eq:kl_property}) is the Kullback-Leibler property and (\ref{eq:Q_star})
has been referred to as the $Q^*$ property by \ctn{Walker04a}.
Assume further that
\begin{equation}
\underset{n}{\sup}~Var\left(\frac{I_{n+1}}{I_n}\right)<\infty.
\label{eq:finite_sup}
\end{equation}
Then,
\begin{equation}
n^{-1}\log\left(I_n\right)\rightarrow -\delta,
\label{eq:bf_consistency_iid}
\end{equation}
almost surely.
\end{theorem}
The following corollary provides the result on asymptotic comparison between two models using
Bayes factors, in the $iid$ case.
\begin{corollary}[\ctn{Walker04a}]
\label{corollary:walker}
Let $R_n(f)=\prod_{i=1}^n\frac{f(X_i)}{f_0(X_i)}$. For $j=1,2$, let
\[
I_{jn}=\int R_n(f)\pi_j(df),
\]
where $\pi_1$ and $\pi_2$ are two different priors on $f$. Let $B_n=I_{1n}/I_{2n}$ denote
the Bayes factor for comparing the two models associated with $\pi_1$ and $\pi_2$. If $\pi_1$ and $\pi_2$  
have the Kullback-Leibler property (\ref{eq:kl_property}) 
with $\delta=\delta_1$ and $\delta=\delta_2$ respectively, satisfy the $Q^*$ property (\ref{eq:Q_star}), and
(\ref{eq:finite_sup}) with $I_n=I_{jn}$, for $j=1,2$, then
\[
n^{-1}\log B_n\rightarrow \delta_2-\delta_1,
\]
almost surely.
\end{corollary}

\begin{remark}
In \ctn{Walker04a} the densities are assumed to be dominated by the Lebesgue measure. However,
this is not necessary. The results remain true if the densities are with respect to any valid
measure; see, for example, \ctn{Barron99} for related concepts and results (Lemma 4 in particular) with respect to general
measures. As such, in our $SDE$-based situation, although the densities are not dominated by the Lebesgue measure
(see (\ref{eq:densities})), all our results still remain valid.
\end{remark}

\subsection{Verification of Theorem \ref{theorem:walker} in $iid$ $SDE$ set-up}
\label{subsec:verification_walker}

In our parametric case, $f_0\equiv f_{\btheta_0}$ and $f\equiv f_{\btheta_1}$.
In this $iid$ set-up, as mentioned earlier $\bxi_j=\xi_{0j}$ for $j=0,1$, so that
$\phi_{\bxi_j}\equiv\xi_{0j}$. For our convenience, we let, for $j=0,1$ and $ i=1,\ldots,n$,
\begin{equation}
\tilde U_{i,\bbeta_j} =\int_0^{T_i}\frac{b_{\bbeta_j}(X_i(s))}{\sigma^2(X_i(s))}dX_i(s),
\quad\quad \tilde V_{i,\bbeta_0,\bbeta_j} =\int_0^{T_i}\frac{b_{\bbeta_0}(X_i(s))b_{\bbeta_j}(X_i(s))}{\sigma^2(X_i(s))}ds.
\label{eq:u_v_2}
\end{equation}
Note that, for $i=1,\ldots,n$, $\tilde V_{i,\bbeta_0}=\tilde V_{i,\bbeta_0,\bbeta_0}$ 
and $\tilde V_{i,\bbeta_1}=\tilde V_{i,\bbeta_1,\bbeta_1}$.
The Kullback-Leibler divergence measure between $f_0$ and $f$ in this set-up is given, with $i=1$, by 
\begin{equation}
\mathcal K(f_{\btheta_0},f_{\btheta_1})=
\frac{\xi_{00}^2}{2}E_{\btheta_0}(\tilde V_{1,\bbeta_0})
-\xi_{00}\xi_{01}E_{\btheta_0}(\tilde V_{1,\bbeta_0,\bbeta_1})
+\frac{\xi_{01}^2}{2}E_{\btheta_0}(\tilde V_{1,\bbeta_1}),
\label{eq:kl1}
\end{equation}
where $E_{\btheta_0}\equiv E_{f_{\btheta_0}}$.
The result easily follows from (\ref{eq:u_v_relation}) and (\ref{eq:zero_mean}).
Now let 
\begin{align}
\delta&=\underset{\bTheta}{\min}~\mathcal K\left(f_{\btheta_0},f_{\btheta_1}\right)\notag\\
&=\underset{\bTheta}{\min}~\left\{\frac{\xi_{00}^2}{2}E_{\btheta_0}(\tilde V_{1,\bbeta_0})
-\xi_{00}\xi_{01}E_{\btheta_0}(\tilde V_{1,\bbeta_0,\bbeta_1})
+\frac{\xi_{01}^2}{2}E_{\btheta_0}(\tilde V_{1,\bbeta_1})\right\}.
\label{eq:delta}
\end{align}
Since $E_{\btheta_0}(\tilde V_{1,\bbeta_1})$ and $E_{\btheta_0}(\tilde V_{1,\bbeta_0,\bbeta_1})$ are 
continuous in $\bbeta_1$, compactness of $\bTheta$ guarantees that
$0\leq \delta<\infty$.

\subsubsection{Verification of (\ref{eq:kl_property})}
\label{subsubsec:kl_property}

To see that (\ref{eq:kl_property}) holds in our case for any prior dominated by Lebesgue measure, 
first let us define
\begin{equation}
\mathcal K^*\left(f_{\overline\btheta},f_{\btheta_1}\right)=\mathcal K\left(f_{\btheta_0},f_{\btheta_1}\right)
-\mathcal K\left(f_{\btheta_0},f_{\overline\btheta}\right),
\label{eq:k_star}
\end{equation}
where $f_{\overline\btheta}=\underset{\bTheta}{\arg\min}~\mathcal K\left(f_{\btheta_0},f_{\btheta_1}\right)$.
Now, let us choose any prior $\pi$ such that $\frac{d\pi}{d\nu}=\varrho$ where $\varrho$ is a continuous positive density
with respect to Lebesgue measure, where, by ``positive" density, we mean a density excluding any interval of null measure.
For any $c_1>0$, we then need to show that
$$\pi\left(\btheta_1\in\bTheta:\delta\leq \mathcal K(f_{\btheta_0},f_{\btheta_1})<\delta+c_1\right)>0,$$ for any prior $\pi$ 
dominated by Lebesgue measure.
This is equivalent to showing 
$$\pi\left(\btheta_1\in\bTheta:0\leq \mathcal K^*(f_{\overline\btheta},f_{\btheta_1})<c_1\right)>0,$$ for any prior $\pi$ 
dominated by Lebesgue measure.

Since $\mathcal K(f_{\btheta_0},f_{\btheta_1})$ is continuous in $\btheta_1$, so is 
$\mathcal K^*(f_{\overline\btheta},f_{\btheta_1})$. 
Compactness of $\bTheta$ ensures uniform continuity of $\mathcal K^*(f_{\overline\btheta},f_{\btheta_1})$.
Hence, for any $c_1>0$, there exists $\epsilon_{c_1}$ independent of $\btheta_1$, such that 
$\|\btheta_1-\overline\btheta\|<\epsilon_{c_1}$ implies $\mathcal K^*(f_{\overline\btheta},f_{\btheta_1})<c_1$.
Then,
\begin{align}
&\pi\left(\btheta_1\in\bTheta:0\leq \mathcal K^*(f_{\overline\btheta},f_{\btheta_1})<c_1\right)
\geq \pi\left(\btheta_1\in\bTheta: \|\btheta_1-\overline\btheta\|<\epsilon_{c_1}\right)\notag\\
&\geq \left[\underset{\left\{\btheta_1\in\bTheta:\|\btheta_1-\overline\btheta\|<\epsilon_{c_1}\right\}}{\inf}~\varrho(\btheta_1)\right]\times
\nu\left(\left\{\btheta_1\in\bTheta:\|\btheta_1-\overline\btheta\|<\epsilon_{c_1}\right\}\right)>0,
\label{eq:positivity1}
\end{align}
where $\nu$ stands for Lebesgue measure.
In other words, (\ref{eq:kl_property}) holds in our case.

\subsubsection{Verification of (\ref{eq:Q_star})}
\label{subsubsec:Q_star}

To see that (\ref{eq:Q_star}) also holds in our $SDE$ set-up, first note that in our case
\begin{equation}
\hat f_{nA(\epsilon)}(x)=\frac{\int_{A(\epsilon)}f_{\btheta_1}(x)\pi_n(d\btheta_1)}
{\int_{A(\epsilon)}\pi_n(d\btheta_1)},
\label{eq:post_pred}
\end{equation}
with 
\begin{align}
A(\epsilon)&=\left\{\btheta_1\in\bTheta:\mathcal K(f_{\btheta_0},f_{\btheta_1})\geq \epsilon\right\}\notag\\
&=\left\{\btheta_1\in\bTheta:
\frac{\xi_{00}^2}{2}E_{\btheta_0}(\tilde V_{1,\bbeta_0})
-\xi_{00}\xi_{01}E_{\btheta_0}(\tilde V_{1,\bbeta_0,\bbeta_1})
+\frac{\xi_{01}^2}{2}E_{\btheta_0}(\tilde V_{1,\bbeta_1})\geq\epsilon\right\}
\label{eq:A_e}
\end{align}
for any $\epsilon>0$. Note that, here we have replaced $K(f_{\btheta_0},f_{\btheta_1})> \epsilon$
with $K(f_{\btheta_0},f_{\btheta_1})\geq \epsilon$ in the definition of $A(\epsilon)$ 
because of continuity of the posterior of $\btheta_1$.
Note that
\begin{equation}
\hat f_{nA(\epsilon)}(X)\leq \underset{\btheta_1\in A(\epsilon)}{\sup}~f_{\btheta_1}(X)
=f_{\hat\btheta_1(X)}(X),
\end{equation}
where $\hat\btheta_1(X)$, which depends upon $X$, is the maximizer lying in the compact set $A(\epsilon)$.
Now note that
\begin{align}
\mathcal K(f_{\btheta_0},\hat f_{nA(\epsilon)})&=E_{\btheta_0}\left[\log f_{\btheta_0}(X)\right]
-E_{\btheta_0}\left[\log \hat f_{nA(\epsilon)}(X)\right]\notag\\
&\geq E_{\btheta_0}\left[\log f_{\btheta_0}(X)\right]-E_{\btheta_0}\left[\log f_{\hat\btheta_1(X)}(X)\right]\notag\\
&=E_{\btheta_0}\left(\log\frac{f_{\btheta_0}(X)}{f_{\hat\btheta_1(X)}(X)}\right).
\label{eq:kl2}
\end{align}
To show that $E_{\btheta_0}\left(\log\frac{f_{\btheta_0}(X)}{f_{\hat\btheta_1(X)}(X)}\right)\geq 0$, we first write (\ref{eq:kl2}) as
\begin{equation}
E_{\btheta_0}\left(\log\frac{f_{\btheta_0}(X)}{f_{\hat\btheta_1(X)}(X)}\right)=
E_{\hat\btheta_1(X)|\btheta_0}E_{X|\hat\btheta_1(X),\btheta_0}\left(\log\frac{f_{\btheta_0}(X)}{f_{\hat\btheta_1(X)}(X)}\right).
\label{eq:kl3}
\end{equation}
In (\ref{eq:kl3}), $E_{X|\hat\btheta_1(X),\btheta}\left(\log\frac{f_{\btheta}(X)}{f_{\hat\btheta_1(X)}(X)}\right)$
is the expectation of $\log\frac{f_{\btheta}(X)}{f_{\hat\btheta_1(X)}(X)}$ with respect to the conditional distribution of
$[X|\hat\btheta_1(X),\btheta]$. Assuming that $Y=\hat\btheta_1(X)$ has density $g_{\btheta}(Y)$ with respect to Lebesgue measure, 
with $X\sim f_{\btheta}$, the aforementioned conditional distribution has density $f_{\theta}(X|Y)=\frac{f_{\btheta}(X)}{g_{\btheta}(Y)}\bI_{\left\{\hat\btheta_1(X)\right\}}(Y)$
(see \ctn{Schervish95} for details).
Hence, letting $P_0$ denote the probability measure associated with $f_{\btheta_0}$, we obtain
\begin{align}
&E_{X|\hat\btheta_1(X)=Y,\btheta_0}\left(\log\frac{f_{\btheta_0}(X)}{f_{\left\{\hat\btheta_1(X)\right\}}(X)}\right)\notag\\
&=\int \log\left(\frac{f_{\btheta_0}(X)}{f_{Y}(X)}\right)\frac{f_{\btheta_0}(X)}{g_{\btheta_0}(Y)}\bI_{\left\{\hat\btheta_1(X)=Y\right\}}(X)dP_0
\notag\\
&=\int \log\left[\frac{\left(f_{\btheta_0}(X)/g_{\btheta_0}(Y)\right)\bI_{\left\{\hat\btheta_1(X)=Y\right\}}(X)}{f_{Y}(X)}\right]
\frac{f_{\btheta_0}(X)}{g_{\btheta_0}(Y)}\bI_{\left\{\hat\btheta_1(X)=Y\right\}}(X)dP_0\notag\\
&\qquad\qquad+\int\log\left(g_{\btheta_0}(Y)\right)\frac{f_{\btheta_0}(X)}{g_{\btheta_0}(Y)}\bI_{\left\{\hat\btheta_1(X)=Y\right\}}(X)dP_0\notag\\
&=\int \log\left(\frac{f_{\btheta_0}(X|Y)}{f_{Y}(X)}\right) f_{\btheta_0}(X|Y)dP_0
+\log\left(g_{\btheta_0}(Y)\right).
\label{eq:kl4}
%
\end{align}
The first term of (\ref{eq:kl4}) is the Kullback-Leibler divergence between the two different densities $f_{\btheta_0}(\cdot|Y)$ and $f_{Y}(\cdot)$, and hence, 
$\mathcal K\left(f_{\btheta_0}(\cdot|Y),f_{Y}(\cdot)\right)>0$, almost surely, for all $Y$. Hence, 
\begin{equation}
E_{Y|\theta_0}E_{X|Y,\btheta_0}\left(\frac{f_{\btheta_0}(X|Y)}{f_{Y}(X)}\right)>0.
\label{eq:kl_5}
\end{equation}
Also, by Jensen's inequality,
\begin{equation}
E_{Y|\theta_0}\left[\log\left(g_{\btheta_0}(Y)\right)\right]=-E_{Y|\theta_0}\left[\log\left(\frac{1}{g_{\btheta_0}(Y)}\right)\right]
\geq -\log E_{Y|\theta_0}\left(\frac{1}{g_{\btheta_0}(Y)}\right).
\label{eq:jensen1}
\end{equation}
Now note that $A(\epsilon)=\left\{\btheta_1\in\bTheta:\mathcal K^*\left(f_{\bar\btheta},f_{\btheta_1}\right)\geq\epsilon-\delta\right\}$, where we must have
$\epsilon\geq\delta\geq 0$. Since $\mathcal K^*\left(f_{\bar\btheta},f_{\btheta_1}\right)$ is uniformly continuous on $\bTheta$, for any $\epsilon^*=\epsilon-\delta>0$,
there exists $\eta=\eta(\epsilon^*)>0$ such that $\mathcal K^*\left(f_{\bar\btheta},f_{\btheta_1}\right)>\epsilon^*$ implies $\|\btheta_1-\bar\btheta\|\geq\eta$.
Let $B(\eta)=\left\{\btheta_1\in\bTheta:\|\btheta_1-\bar\btheta\|\geq\eta\right\}$.
It follows that
\begin{equation*}
E_{Y|\theta_0}\left(\frac{1}{g_{\btheta_0}(Y)}\right)=\left|A(\epsilon)\right|=\int_{A(\epsilon)}dy\leq \int_{B(\eta)}dy=\left|B(\eta)\right|.
\end{equation*}
Now we can achieve $\left|B(\eta)\right|<\exp(-\epsilon)$ by suitable 
reparameterization of the components of $\btheta_0$ and $\btheta_1$ lying in the compact space $\bTheta$. For instance, if $\theta_{1j}$, the $j$-th component of $\btheta_1$ 
satisfies $a_j\leq\theta_{1j}\leq b_j$, then for any $c_j\geq 1$,
$\theta_{1j}=c_j\tilde\theta_{1j}$, where $\tilde\theta_{1j}=\theta_{1j}/c_j\in [a_jc^{-1}_j,b_jc^{-1}_j]$. 
We also write $\theta_{0j}=c_j\tilde\theta_{0j}$, where $\tilde\theta_{0j}=\theta_{0j}/c_j$; here $\theta_{0j}$ is the $j$-th component of $\btheta_0$. 
Abusing notation, we continue denote the parameter space associated with the reparameterizations $\tilde\theta_{1j}$ by $\bTheta$. 
By choosing $c_j$'s to be sufficiently large, the inequality $\left|B(\eta)\right|<\exp(-\epsilon)$ 
can be easily achieved.
The interpretation of this is that the part of the parameter space with large Kullback-Leibler divergence from the true density has relatively small volume. 

Then, for both the cases it follows from (\ref{eq:jensen1}) that
$E_{Y|\theta_0}\left[\log\left(g_{\btheta_0}(Y)\right)\right]\geq \epsilon$. Combining this with (\ref{eq:kl_5}), (\ref{eq:kl4}), (\ref{eq:kl3}) and (\ref{eq:kl2}), we have that
$\mathcal K(f_{\btheta_0},\hat f_{nA(\epsilon)})\geq\epsilon$.
Hence, (\ref{eq:Q_star}) is satisfied in our $SDE$ set-up.

\subsubsection{Verification of (\ref{eq:finite_sup})}
\label{subsubsec:finite_sup}

We now prove that (\ref{eq:finite_sup}) also holds. 
It is straightforward to verify that 
\begin{equation}
\frac{I_{n+1}}{I_n}=\frac{\hat f_{n+1}(X_{n+1})}{f_{\btheta_0}(X_{n+1})},
\label{eq:I_ratio}
\end{equation}
where
\begin{equation}
\hat f_{n+1}(\cdot)=E_{\btheta_1|X_1,\ldots,X_n}\left[f_{\btheta_1}(\cdot)\right]
\label{eq:hat_f}
\end{equation}
is the posterior predictive distribution of $f_{\btheta_1}(\cdot)$, with respect to the posterior
of $\btheta_1$, given $X_1,\ldots,X_n$. In (\ref{eq:hat_f}), $E_{\btheta_1|X_1,\ldots,X_n}$ denotes expectation with respect
to the posterior of $\btheta_1$ given $X_1,\ldots,X_n$. 



First note that, since
\begin{equation}
\log\left[f_{\btheta_0}(X_{n+1})\right]=\xi_{00}\tilde U_{n+1,\bbeta_0}-\frac{\xi_{00}^2}{2}\tilde V_{n+1,\bbeta_0},
\label{eq:log_f_0}
\end{equation}
it follows from Lemma \ref{lemma:moment_existence} that the moments of all orders of 
$\log\left[f_{\btheta_0}(X_{n+1})\right]$ exist and are finite. 
Also, since $X_i$ are $iid$, the moments are the same for every $n=1,2,\ldots$.
In other words, 
\begin{equation}
\underset{n}{\sup}~Var\left(\log f_{\btheta_0}(X_{n+1})\right)<\infty.
\label{eq:var_deno_finite}
\end{equation}

Then observe that for any given $X_{n+1}$, using compactness of $\bTheta$ and continuity of $f_{\btheta_1}(X_{n+1})$
with respect to $\btheta_1$, 
$$f_{\btheta^{*}_1(X_{n+1})}(X_{n+1})=\underset{\btheta_1\in\bTheta}{\inf}~f_{\btheta_1}(X_{n+1})\leq
\hat f_{n+1}(X_{n+1})\leq \underset{\btheta_1\in\bTheta}{\sup}~f_{\btheta_1}(X_{n+1})=f_{\btheta^{**}_1(X_{n+1})}(X_{n+1}),$$
where $\btheta^*_1(X_{n+1})= \underset{\btheta_1\in\bTheta}{\arg\min}~f_{\btheta_1}(X_{n+1})$ and
$\btheta^{**}_1(X_{n+1})= \underset{\btheta_1\in\bTheta}{\arg\max}~f_{\btheta_1}(X_{n+1})$. 
Clearly,
$\btheta^*_1(X_{n+1}), \btheta^{**}_1(X_{n+1})\in\bTheta$, for any given $X_{n+1}$. Moreover,
$$\btheta^*_1(X_{n+1})=(\bbeta^*_1(X_{n+1}),\xi^*_{01}(X_{n+1})),~~~~ 
\btheta^{**}_1(X_{n+1})=(\bbeta^{**}_1(X_{n+1}),\xi^{**}_{01}(X_{n+1})))$$ where 
each component of $\btheta^*_1(X_{n+1})$ and $\btheta^{**}_1(X_{n+1})$ depends on $X_{n+1}$.
Noting that $U_{n+1,\btheta^*_1(X_{n+1})}=\xi^*_{01}(X_{n+1})\tilde U_{n+1,\bbeta^*_1(X_{n+1})}$
and $V_{n+1,\btheta^*_1(X_{n+1})}=\left\{\xi^*_{01}(X_{n+1})\right\}^2\tilde V_{n+1,\bbeta^*_1(X_{n+1})}$,
it follows from the above inequality that
\begin{align}
&-\left|U_{n+1,\btheta^*_1(X_{n+1})}\right| 
-\frac{V_{n+1,\btheta^*_1(X_{n+1})}}{2}\notag\\ 
&\leq U_{n+1,\btheta^*_1(X_{n+1})} -\frac{V_{n+1,\btheta^*_1(X_{n+1})}}{2}\notag\\
&\leq\log\hat f_{n+1}(X_{n+1})\notag\\
&\leq U_{n+1,\btheta^{**}_1(X_{n+1})}
-\frac{V_{n+1,\btheta^{**}_1(X_{n+1})}}{2}\notag\\
&\leq \left|U_{n+1,\btheta^{**}_1(X_{n+1})}\right|
+\frac{V_{n+1,\btheta^{**}_1(X_{n+1})}}{2}.\notag 
\end{align}
Hence, $E_{\btheta_0}\left(\log\hat f_{n+1}(X_{n+1})\right)^2$ lies between
$E_{\btheta_0}\left(\left|U_{n+1,\btheta^*_1(X_{n+1})}\right|
+\frac{V_{n+1,\btheta^*_1(X_{n+1})}}{2}\right)^2$ and
\\
$E_{\btheta_0}\left(\left|U_{n+1,\btheta^{**}_1(X_{n+1})}\right|
+\frac{V_{n+1,\btheta^{**}_1(X_{n+1})}}{2}\right)^2$. 

We obtain uniform lower and upper bounds of the above two expressions
in the following manner. For the upper bound of the latter
%
%
we first take supremum of the expectation
with respect to $X_{n+1}$, conditional on $\btheta^{**}_1(X_{n+1})=\bvarsigma$, over $\bvarsigma\in\bTheta$, and then take 
expectation with respect to
$X_{n+1}$. 
Since $\bvarsigma\in\bTheta$, compactness of $\bTheta$ and Lemma \ref{lemma:moment_existence}
ensure that  
the moments of any given order of the above expression is uniformly bounded above. 
Analogously, we obtain a uniform lower bound replacing the supremum with infimum. 
In the same way we obtain uniform lower and upper bounds of the other expression.
The uniform bounds on the second order moments, in turn, guarantee that
\begin{equation}
\underset{n}{\sup}~Var\left(\log \hat f_{n+1}(X_{n+1})\right)<\infty.
\label{eq:var_num_finite}
\end{equation}
Combining (\ref{eq:var_deno_finite}) and (\ref{eq:var_num_finite}) and using the Cauchy-Schwartz inequality
for the covariance term associated with $Var\left(\log \hat f_{n+1}(X_{n+1})-\log f_{\btheta_0}(X_{n+1})\right)$
shows that (\ref{eq:finite_sup}) holds in our set-up.

We formalize the above arguments in the form of a theorem in the $SDE$ based $iid$ set-up. 
\begin{theorem}
\label{theorem:bf_consistency_iid}
Assume the $iid$ case of the $SDE$ based set-up and conditions (H1$^\prime$) -- (H4$^\prime$).
Then (\ref{eq:bf_consistency_iid}) holds.
\end{theorem}
The following corollary in the $iid$ $SDE$ context is motivated by Corollary \ref{corollary:walker}.
\begin{corollary}
\label{corollary:sde_iid}
For $j=1,2$, let $R_{jn}(\btheta_j)=\prod_{i=1}^n\frac{f_{\btheta_j}(X_i)}{f_{\btheta_0}(X_i)}$, where
$\btheta_1$ and $\btheta_2$ are two different finite sets of parameters,
perhaps with different dimensionalities,
associated with the two models
to be compared. For $j=1,2$, let
\[
I_{jn}=\int R_{jn}(\btheta_j)\pi_j(d\btheta_j),
\]
where $\pi_j$ is the prior on $\btheta_j$. 
Let $B_n=I_{1n}/I_{2n}$ as before. 
Assume the $iid$ case of the $SDE$ based set-up and suppose that both
the models satisfy conditions (H1$^\prime$) -- (H4$^\prime$) and
have the Kullback-Leibler property 
with $\delta=\delta_1$ and $\delta=\delta_2$ respectively. 
Then 
\[
n^{-1}\log B_n\rightarrow \delta_2-\delta_1,
\]
almost surely.
\end{corollary}

\section{General asymptotic theory of Bayes factor in the non-$iid$ set-up}
\label{sec:non_iid}

In this section, we first develop a general asymptotic theory of Bayes factors in the non-$iid$ set-up, and then
obtain the result for the non-$iid$ $SDE$ set-up as a special case of our general theory.

\subsection{The basic set-up}
\label{eq:non_iid_set-up}


We assume that for $i=1,\ldots,n$, $X_i\sim f_ {0i}$, that is, the true density function corresponding to 
the $i$-th individual is $f_{0i}$. Considering another arbitrary density $f_i$ for individual $X_i$ we 
investigate consistency of the 
Bayes factor in this general non-$iid$ set-up.
For our purpose we introduce the following two properties:
\\[2mm]
{\bf 1.  Kullback-Leibler $(\delta)$ property in the non-$iid$ set-up}:
\\[2mm]
We denote the Kullback-Leibler divergence measure between $f_{0i}$ and $f_i$ by $\mathcal K(f_{0i},f_i)$ 
and assume that the limit
\begin{equation}
\mathcal K^{\infty}\left(f_0,f\right)
=\underset{n\rightarrow\infty}{\lim}~\frac {1}{n} \sum_{i=1}^n E\left[\log\frac{f_{0i}(X_i)}{f_i(X_i)}\right]
=\underset{n\rightarrow\infty}{\lim}~\frac {1}{n} \sum_{i=1}^n\mathcal K\left(f_{0i},f_i\right)
\label{eq:kl_infinite}
\end{equation}
exists almost surely with respect to the prior $\pi$ on $f$.
Let the prior distribution $\pi$ satisfy
\begin{equation}
\pi\left(f:\underset{i}{\inf}~\mathcal K(f_{0i},f_i)\geq\delta\right)=1,
\label{eq:prior_delta}
\end{equation}
for some $\delta\geq 0$.
Then we say that $\pi$ has the Kullback-Leibler $(\delta)$ property if, for any $c>0$, 
\begin{equation}
\pi\left(f:\delta\leq\mathcal K^{\infty}\left(f_0,f\right)\leq\delta+c\right)>0.
\label{eq:kl_property_non_iid}
\end{equation}
{\bf 2. $Q^*$ property in the non-$iid$ set-up}:
\\[2mm]
Let us denote the posterior distribution corresponding to $n$ observations by $\pi_n$. 
We denote $\pi(df_1,df_2,\ldots,df_n)$ by $\pi(\tilde {df})$. For any set $A$,
\[
\pi_n(A)=\frac{\int_A\prod_{i=1}^nf_i(X_i)\pi(\tilde {df})}{\int\prod_{i=1}^nf_i(X_i)\pi(\tilde {df})}
\]
denotes the posterior probability of $A$.
Let $$R_n(f_1,f_2,\ldots,f_n)=\prod_{i=1}^n \frac{f_i(X_i)}{f_{0i}(X_i)}.$$
Let us define the posterior predictive density by
\[
\hat f_n(X_n)=\int f_n(X_n)\pi_n({df}_n),
\]
and
\[
\hat f_{nA} (X_n)=\int f_n(X_n)\pi_{nA}({df}_n)
\]
to be the posterior predictive density with posterior restricted to the set $A$, that is, for $\pi(A)>0$,
\[
\pi_{nA}({df}_n)=\frac{{\bI}_A(f_n)\pi_n({df}_n)}{\int_A \pi_n({df}_n)}.
\]
Then we say that the prior has the {\it property $Q^*$ in the non-$iid$ set-up} if the following holds for any $\epsilon>0$:
\begin{equation}
\underset{n}{\lim\inf}~ \mathcal K(f_{0n},\hat f_{n,A_n(\epsilon)})\geq \epsilon,
\label{eq:Q_star_non_iid}
\end{equation}
when 
\begin{equation}
A_n(\epsilon)=\{f_n:\mathcal K(f_{0n},f_n)\geq \epsilon\}.
\label{eq:A_n}
\end{equation}
Let $I_0\equiv 1$ and for $n\geq 1$, let us define
\begin{equation}
I_n=\int R_n(f_1,f_2,\ldots,f_n)\pi(\tilde{df}), 
\label{eq:bf_non_iid}
\end{equation}
which is relevant for the study of the Bayes factors. 
Regarding convergence of $I_n$, we formulate the following theorem.
\begin{theorem}
\label{theorem:bf_non_iid}
Assume the non-$iid$ set-up and that the limit (\ref{eq:kl_infinite}) exists almost surely with respect to the prior $\pi$.
Also assume that the prior $\pi$ satisfies (\ref{eq:prior_delta}), has the Kullback-Leibler $(\delta)$ and $Q^*$ properties
given by (\ref{eq:kl_property_non_iid}) and (\ref{eq:Q_star_non_iid}), respectively. 
Assume further that
\begin{equation}
\underset{i}{\sup}~ E\left[\log\frac{f_{0i}(X_i)}{f_i(X_i)}\right]^2<\infty 
\label{eq:sup1}
\end{equation}
and
\begin{equation}
\underset{n}{\sup}~E\left[\log\frac{I_n}{I_{n-1}}\right]^2<\infty.
\label{eq:convergent_series2}
\end{equation}
Then 
\begin{equation}
n^{-1}\log I_n \rightarrow -\delta,
\label{eq:bf_consistency_non_iid}
\end{equation}
almost surely as $n\rightarrow\infty$.
\end{theorem}

\begin{corollary}
\label{corollary:bf_non_iid}
For $j=1,2$, let
\[
I_{jn}=\int R_n(f_1,\ldots,f_n)\pi_j(\tilde df),
\]
where $\pi_1$ and $\pi_2$ are two different priors on $f$. 
Let $B_n=I_{1n}/I_{2n}$ denote
the Bayes factor for comparing the two models associated with $\pi_1$ and $\pi_2$. If both the models satisfy
the conditions of Theorem \ref{theorem:bf_non_iid}, and satisfy the Kullback-Leibler property  
with $\delta=\delta_1$ and $\delta=\delta_2$ respectively, then
\[
n^{-1}\log B_n\rightarrow \delta_2-\delta_1,
\]
almost surely.
\end{corollary}

\section{Specialization of non-$iid$ asymptotic theory of Bayes factors to non-$iid$ $SDE$ set-up
where $T_i$ are constants for every $i$ but $n\rightarrow\infty$}
\label{sec:sde_non_iid}

In this section we relax the restrictions $T_i=T$
and $x^i=x$ for $i=1,\ldots,n$. In other words, here we deal with the set-up
where the processes $X_i(\cdot);~i=1,\ldots,n$, are independently,
but not identically distributed.
Following \ctn{Maitra14a}, \ctn{Maitra14b} we assume the following:
\begin{itemize}
\item[(H6$^\prime$)] The sequences $\{T_1,T_2,\ldots\}$ and 
$\{x^1,x^2,\ldots\}$ are sequences in compact sets $\mathfrak T$ and $\mathfrak X$, respectively, so that
there exist convergent subsequences with limits in $\mathfrak T$ and $\mathfrak X$.
For notational convenience, we continue to denote the convergent subsequences as  $\{T_1,T_2,\ldots\}$
and $\{x^1,x^2,\ldots\}$. Let us denote the limits by $T^{\infty}$ and $x^{\infty}$, where $T^{\infty}\in\mathfrak T$
and $x^{\infty}\in\mathfrak X$. 
\end{itemize}
\begin{remark}
Note that the choices of the convergent subsequences $\{T_1,T_2,\ldots\}$ and 
$\{x^1,x^2,\ldots\}$ are not unique. However, this non-uniqueness does not affect asymptotic selection
of the correct model via Bayes factor. Indeed, as will be evident from our proof, 
for any choice of convergent subsequence, the Bayes factor almost surely converges exponentially 
to the correct quantity.
The reason for this is that we actually need to deal with the infimum of the Kullback-Leibler distance
over $\mathfrak X$ and $\mathfrak T$, which is of course independent of the choices of subsequences; see
Section \ref{subsec:kl_property_non_iid} for the details.
\end{remark}

Following \ctn{Maitra14a}, we denote the process associated with the initial value $x$ and time point $t$ as $X(t,x)$, so that
$X(t,x^i)=X_i(t)$, and $X_i=\left\{X_i(t);~t\in[0,T_i]\right\}$.

Let $\btheta_j=(\bbeta_j,\bxi_j)$ for $j=0,1$ denote the set of finite number of parameters, where
$\bbeta_j$ and $\bxi_j$ have the same interpretation as in the $iid$ set-up. 
As before, $\bz_i(t)=(z_{i1}(t),z_{i2}(t),\ldots,z_{ip}(t))$ is the set of covariate information 
corresponding to $i$-th individual
at time point $t$.
For $x^i\in \mathfrak X$, $T_i\in \mathfrak T$, $\bz_i(t)\in\boldcal Z$ and $\btheta_j\in\bTheta$, let
\begin{align}
U_{x^i,T_i,\bz_i,\btheta_j}
&=\int_0^{T_i}\frac{\phi_{i,\bxi_j}(s)b_{\bbeta_j}(X_i(s,x^i))}{\sigma^2(X_i(s,x^i))}d X_i(s,x^i)
\label{eq:u_x_T_non_iid};\\
V_{x^i,T_i,\bz_i,\btheta_0,\btheta_j}&=\int_0^{T_i}\frac{\phi_{i,\bxi_j}(s)\phi_{i,\bxi_0}(s)b_{\bbeta_j}(X_i(s,x^i))
b_{\bbeta_0}(X_i(s,x^i))}{\sigma^2(X_i(s,x^i))}ds.
\label{eq:v_x_T_non_iid}
\end{align}
As before, 
$V_{x^i,T_i,\bz_i,\btheta_0}=V_{x^i,T_i,\bz_i,\btheta_0,\btheta_0}$ 
and $V_{x^i,T_i,\bz_i,\btheta_1}=V_{x^i,T_i,\bz_i,\btheta_1,\btheta_1}$.


In this non-$iid$ set-up $f_{0i}=f_{\btheta_0,x^i,T_i,\bz_i}$ and $f_{i}=f_{\btheta_1,x^i,T_i,\bz_i}$. 
An extension of Lemma \ref{lemma:moment_existence} incorporating $x$, $T$ and $\bz$ shows that 
moments of $U_{x,T,\bz,\btheta_j}$, $V_{x,T,\bz,\btheta_j}$, 
$V_{x,T,\bz,\btheta_0,\btheta_j}$ of all orders exist, and are 
continuous in $x$, $T$, $\bz$, $\btheta_1$. Formally, we have the following lemma. 
\begin{lemma}
\label{lemma:moment_existence_non_iid}
Assume (H1$^\prime$) -- (H6$^\prime$). Then for all $x\in\mathfrak X$, $T\in\mathfrak T$, $\bz\in\boldcal Z$ 
and $\btheta_1\in\bTheta$, for $k\geq 1$,  
\begin{align}
E_{\btheta_0}\left[U_{x,T,\bz,\btheta_j}\right]^k&<\infty; j=0,1,\label{eq:finite_moment_u_non_iid}\\
E_{\btheta_0}\left[V_{x,T,\bz,\btheta_1}\right]^k&<\infty,\label{eq:finite_moment_v1_non_iid}\\
E_{\btheta_0}\left[V_{x,T,\bz,\btheta_0,\btheta_j}\right]^k&<\infty; j=0,1.\label{eq:finite_moment_v_non_iid}
\end{align}
Moreover, the above expectations are continuous in $(x,T,\bz,\btheta_1)$.
\end{lemma}

In particular, 
the Kullback-Leibler distance
is continuous in $x$, $T$, $\bz$ and $\btheta_1$.
The following lemma asserts that the average of the Kullback-Leibler distance is also a Kullback-Leibler
distance in the limit.
\begin{lemma}
\label{lemma:kl_average}
The limiting average
$\underset{n\rightarrow\infty}{\lim}~\frac{1}{n}\sum_{k=1}^n\mathcal K(f_{\btheta_0,x^k,T_k,\bz_k},
f_{\btheta_1,x^k,T_k,\bz_k})$
is also a Kullback-Leibler distance.
\end{lemma}

Even in this non-$iid$ context, the Bayes factor is of the same form as (\ref{eq:I_n}); however,
for $j=0,1$, $U_{x^i,T_i,\bz_i,\bbeta_j,\bxi_j}$ and $V_{x^i,T_i,\bz_i,\bbeta_j,\bxi_j}$ 
are not identically distributed for $i=1,\ldots,n$.
Next, we establish strong consistency of Bayes factor in the non-$iid$ $SDE$ set-up by verifying the sufficient conditions
of Theorem \ref{theorem:bf_non_iid}.

\subsection{Verification of (\ref{eq:prior_delta}) and the Kullback-Leibler property in the non-$iid$ set-up}
\label{subsec:kl_property_non_iid}

Firstly, note that in our case, 
\begin{equation}
\mathcal K^{\infty}\left(f_0,f\right)
=\mathcal K^{\infty}\left(f_{\btheta_0},f_{\btheta_1}\right),
\label{eq:kl_limit_2}
\end{equation}
where the rightmost side, as asserted by Lemma \ref{lemma:kl_average}, 
clearly exists almost surely with respect to $\btheta_1$
and is also continuous in $\btheta_1$.

Now note that compactness of $\mathfrak X$, $\mathfrak T$ and  $\boldcal Z$ 
along with continuity of the function $\phi_{\bxi_j}$ and $\mathcal K(f_{\btheta_0,x,T,\bz},f_{\btheta_1,x,T,\bz})$ 
with respect to $x, T$ and $\bz$ implies
\begin{align}
\psi(\btheta_1)&=\underset{x\in\mathfrak X,~T\in\mathfrak T,~\bz\in\boldcal Z}{\inf}~
\mathcal K(f_{\btheta_0,x,T,\bz},f_{\btheta_1,x,T,\bz})\notag\\
&=\underset{x\in\mathfrak X,~T\in\mathfrak T,~\bz\in\boldcal Z}{\inf}~\int_0^T\left\{\frac{\phi^2_{\bxi_0}(\bz(s))}{2}
E_{\btheta_0}(\breve V_{x,\bbeta_0}(s))\right.\notag\\
&\left.\qquad\qquad-\phi_{\bxi_0}(\bz(s))\phi_{\bxi_1}(\bz(s))E_{\btheta_0}(\breve V_{x,\bbeta_0,\bbeta_1}(s))
+\frac{\phi^2_{\bxi_1}(\bz(s))}{2}
E_{\btheta_0}(\breve V_{x,\bbeta_1}(s))\right\}ds\notag\\
&=\underset{x\in\mathfrak X,~T\in\mathfrak T,~\bz(\tilde s(T))\in\boldcal Z}{\inf}
~T\left\{\frac{\phi^2_{\bxi_0}(\bz(\tilde s(T)))}{2}
E_{\btheta_0}(\breve V_{x,\bbeta_0}(\tilde s(T)))\right.\notag\\
&\left.\qquad\qquad-\phi_{\bxi_0}(\bz(\tilde s(T)))\phi_{\bxi_1}(\bz(\tilde s(T)))
E_{\btheta_0}(\breve V_{x,\bbeta_0,\bbeta_1}(\tilde s(T)))
+\frac{\phi^2_{\bxi_1}(\bz(\tilde s(T)))}{2}
E_{\btheta_0}(\breve V_{x,\bbeta_1}(\tilde s(T)))\right\},
\label{eq:kl5}
\end{align}
by the mean value theorem for integrals, where $\tilde s(T)\in[0,T]$
such that the above equality holds.
Also note that the expression in (\ref{eq:kl5}) is continuous in $T$ since originally
the integral on $[0,T]$ is continuous in $T$. 
Now note that if $|T-\tilde T|<\delta_1 (\epsilon)$ such that
$|\phi_{\bxi_j}(\bz(\tilde s(T)))-\phi_{\bxi_j}(\bz(\tilde s(\tilde T)))|<\frac{\epsilon}{2}$
due to continuity in $T$ and if $|\bz(\tilde s(\tilde T))-\tilde \bz(\tilde s(\tilde T))|<\delta_2 (\epsilon)$
such that $|\phi_{\bxi_j}(\bz(\tilde s(\tilde T)))-\phi_{\bxi_j}(\tilde \bz(\tilde s(\tilde T)))|<\frac{\epsilon}{2}$
due to continuity of $\phi_{\bxi_j}$ in $\bz$, then
$|\phi_{\bxi_j}(\bz(\tilde s(T)))-\phi_{\bxi_j}(\tilde\bz(\tilde s(\tilde T)))|
\leq |\phi_{\bxi_j}(\bz(\tilde s(T)))-\phi_{\bxi_j}(\bz(\tilde s(\tilde T)))| 
+|\phi_{\bxi_j}(\bz(\tilde s(\tilde T)))-\phi_{\bxi_j}(\tilde \bz(\tilde s(\tilde T)))|<\epsilon$, showing
that $\phi_{\bxi_j}(\bz(\tilde s(T)))$ is continuous in $T$ and $\bz(\tilde s(T))$, which also belong to compact spaces.
Hence, from (\ref{eq:kl5}) it follows that
\begin{align}
\psi(\btheta_1)&=T^*(\btheta_1)\left[\frac{\phi^2_{\bxi_0}(\bz^*(\btheta_1)(\tilde s(T^*(\btheta_1))))}{2}
E_{\btheta_0}(\breve V_{x^{*}(\btheta_1),\bbeta_0}(\tilde s(T^*(\btheta_1))))\right.\notag\\
&\qquad\qquad\left. -\phi_{\bxi_0}(\bz^*(\btheta_1)(\tilde s(T^*(\btheta_1))))\phi_{\bxi_1}(\bz^*(\btheta_1)(\tilde s(T^*(\btheta_1))))
E_{\btheta_0}(\breve V_{x^{*}(\btheta_1),\bbeta_0,\bbeta_1}(\tilde s(T^*(\btheta_1))))\right.\notag\\
&\qquad\qquad\qquad \left. +\frac{\phi^2_{\bxi_1}(\bz^*(\btheta_1)(\tilde s(T^*(\btheta_1))))}{2}
 E_{\btheta_0}(\breve V_{x^{*}(\btheta_1),\bbeta_1}(\tilde s(T^*(\btheta_1))))\right],
\label{eq:kl_non_iid2}
\end{align}
where $x^{*}(\btheta_1)\in\mathfrak X$, $T^*(\btheta_1)\in\mathfrak T$, 
$\bz^*(\btheta_1)(\tilde s(T^*(\btheta_1)))\in\boldcal Z$ depend upon $\btheta_1$. Then, 
considering the constant correspondence function
$\bgamma(\btheta_1)=\mathfrak X\times\mathfrak T\times\boldcal Z$, for all $\btheta_1\in\bTheta$,
we note that $\bgamma$ is both upper and lower hemicontinuous
(hence continuous), and also compact-valued. Hence, Berge's maximum theorem (\ctn{Berge63}) guarantees that 
(\ref{eq:kl_non_iid2}) is a continuous function of $\btheta_1$. 

Because of continuity of $\psi(\btheta_1)$ in $\btheta_1$, the set $\left\{\btheta_1:\psi(\btheta_1)\geq\delta\right\}$
is open and can be assigned any desired probability by choosing appropriate priors dominated by the Lebesgue measure.
That is, we can assign prior probability one to this set by choosing appropriate priors dominated by the Lebesgue measure.
Now, because of the inequality
\begin{equation*}
\pi\left(\btheta_1: \underset{i}{\inf}~
\mathcal K\left(f_{\btheta_0,x^i,T_i, \bz_i},f_{\btheta_1,x^i,T_i, \bz_i}\right)\geq\delta\right)
\geq \pi\left(\btheta_1:\psi(\btheta_1)\geq\delta\right),
\end{equation*}
and since we choose $\pi$ such that $\pi\left(\btheta_1:\psi(\btheta_1)\geq\delta\right)=1$, it follows that
$$\pi\left(\btheta_1: \underset{i}{\inf}~
\mathcal K\left(f_{\btheta_0,x^i,T_i, \bz_i},f_{\btheta_1,x^i,T_i, \bz_i}\right)\geq\delta\right)=1,$$ 
satisfying (\ref{eq:prior_delta}).

The Kullback-Leibler property of the Lebesgue measure dominated $\pi$ easily follows from continuity of 
$\mathcal K^{\infty}\left(f_{\btheta_0},f_{\btheta_1}\right)$
in $\btheta_1$.

\subsection{Verification of the $Q^*$ property in the non-$iid$ set-up}
\label{subsec:Q_star}

Observe that in this situation, for any $\epsilon>0$,
\begin{align}
A_n(\epsilon)&=\left\{f_n:\mathcal K\left(f_{0n},f_n\right)\geq\epsilon\right\}\notag\\
&=\left\{\btheta_1:\mathcal K\left(f_{\btheta_0,x^n,T_n, \bz_n},f_{\btheta_1,x^n,T_n, \bz_n}\right)\geq\epsilon\right\}\notag
\end{align}
Then note that
\begin{equation}
\hat f_{nA_n(\epsilon)}(X)\leq \underset{\btheta_1\in A_n(\epsilon)}{\sup}~f_{\btheta_1,x^n,T_n, \bz_n}(X)
=f_{\hat\btheta_1(X,x^n,T_n,\bz_n)}(X),
\end{equation}
where $\hat\btheta_1(X,x^n,T_n,\bz_n)$, which depends upon $X,x^n,T_n,\bz_n$, 
is the maximizer lying in the compact set $A_n(\epsilon)$.
Now, 
\begin{align}
\mathcal K(f_{\btheta_0,x^n,T_n,\bz_n},\hat f_{nA_n(\epsilon)})
&=E_{\btheta_0}\left[\log f_{\btheta_0,x^n,T_n,\bz_n}(X)\right]
-E_{\btheta_0}\left[\log \hat f_{nA_n(\epsilon)}(X)\right]\notag\\
&\geq E_{\btheta_0}\left[\log f_{\btheta_0,x^n,T_n,\bz_n}(X)\right]
-E_{\btheta_0}\left[\log f_{\hat\btheta_1(X,x^n,T_n,\bz_n)}(X)\right]\notag\\
&=E_{\btheta_0}\left(\log\frac{f_{\btheta_0,x^n,T_n,\bz_n}(X)}{f_{\hat\btheta_1(X,x^n,T_n,\bz_n)}(X)}\right).
\label{eq:kl_non_iid}
%
\end{align}
In the same way as in Section \ref{subsubsec:Q_star}, after suitable reparameterization, we can achieve $\underset{n}{\sup}~|A_n(\epsilon)|<\exp(-\epsilon)$.
Then as before it can be shown that (\ref{eq:kl_non_iid}) is at least $\epsilon$. 
Hence, (\ref{eq:Q_star_non_iid}) is satisfied in our non-$iid$ $SDE$ set-up.

\subsection{Verification of (\ref{eq:sup1})}
\label{subsec:finite_series_1_non_iid}

From Lemma \ref{lemma:moment_existence_non_iid} it follows that 
$E\left\{\log\frac{f_{\btheta_0,x,T,\bz}(X)}{f_{\btheta_1,x,T,\bz}(X)}\right\}^2$
exists and is continuous in $\btheta_1$, $x$, $T$ and $\bz$. Then compactness of $\bTheta$, 
$\mathfrak X$, $\mathfrak T$ and $\boldcal Z$ ensures
(\ref{eq:sup1}).

\subsection{Verification of (\ref{eq:convergent_series2})}
\label{subsec:finite_series_2_non_iid}

For the non-$iid$ case, the following identity holds:
\begin{align}
\frac{I_{n+1}}{I_n}&=\frac{\hat f_{n+1}(X_{n+1})}{f_{0,n+1}(X_{n+1})}\notag\\
&=\frac{\hat f_{x^{n+1},T_{n+1},\bz_{n+1}}(X_{n+1})}{f_{\btheta_0,x^{n+1},T_{n+1},\bz_{n+1}}(X_{n+1})},
\label{eq:I_ratio_non_iid}
\end{align}
where
\begin{equation}
\hat f_{x^{n+1},T_{n+1}, \bz_{n+1}}(\cdot)=E_{\btheta_1|X_1,\ldots,X_n}\left[f_{\btheta_1,x^{n+1},T_{n+1}, \bz_{n+1}}(\cdot)\right]
\label{eq:hat_f_non_iid}
\end{equation}
is the posterior predictive distribution of $f_{\btheta_1,x^{n+1},T_{n+1},\bz_{n+1}}(\cdot)$, with respect to the posterior
of $\btheta_1$, given $X_1,\ldots,X_n$. 

Now since $\log f_{\btheta_0,x^{n+1},T_{n+1},\bz_{n+1}}(X_{n+1})=
U_{x^{n+1},T_{n+1},\bz_{n+1},\btheta_0}-\frac{V_{x^{n+1},T_{n+1},\bz_{n+1},\btheta_0}}{2}$,
using Lemma \ref{lemma:moment_existence_non_iid} and compactness of $\bTheta$, $\mathfrak X$,
$\mathfrak T$ and $\boldcal Z$ it is easy to see that the moments of $\log f_{\btheta_0,x^{n+1},T_{n+1},\bz_{n+1}}(X_{n+1})$ are
uniformly bounded above. So, we have
\begin{equation}
\underset{n}{\sup}~E\left(\log f_{\btheta_0,x^{n+1},T_{n+1},\bz_{n+1}}(X_{n+1})\right)^2<\infty.
\label{eq:var_deno_finite_non_iid}
\end{equation}

As in the $iid$ case, here also we have
\begin{align}
f_{\btheta^{*}_1(X_{n+1},x^{n+1},T_{n+1},\bz_{n+1})}(X_{n+1})
&=\underset{\btheta_1\in\bTheta}{\inf}~f_{\btheta_1,x^{n+1},T_{n+1},\bz_{n+1}}(X_{n+1})\notag\\
&\leq
\hat f_{x^{n+1},T_{n+1},\bz_{n+1}}(X_{n+1})\notag\\
&\leq \underset{\btheta_1\in\bTheta}{\sup}~f_{\btheta_1,x^{n+1},T_{n+1},\bz_{n+1}}(X_{n+1})
=f_{\btheta^{**}_1(X_{n+1},x^{n+1},T_{n+1},\bz_{n+1})}(X_{n+1}),\notag
\end{align}
where $\btheta^*_1(X_{n+1},x^{n+1},T_{n+1},\bz_{n+1})= 
\underset{\btheta_1\in\bTheta}{\arg\min}~f_{\btheta_1,x^{n+1},T_{n+1},\bz_{n+1}}(X_{n+1})\in\bTheta$ and
$\btheta^{**}_1(X_{n+1},x^{n+1},T_{n+1},\bz_{n+1})= 
\underset{\btheta_1\in\bTheta}{\arg\max}~f_{\btheta_1,x^{n+1},T_{n+1},\bz_{n+1}}(X_{n+1})\in\bTheta$. 
Note that each component of $\btheta^*_1(X_{n+1},x^{n+1},T_{n+1},\bz_{n+1})$ and $\btheta^{**}_1(X_{n+1},x^{n+1},
T_{n+1},\bz_{n+1})$ depends on $X_{n+1},x^{n+1},T_{n+1},\bz_{n+1}$. 

It follows, as in the $iid$ case, that
\begin{align}
&-\left|U_{\btheta^*_1,x^{n+1},T_{n+1},\bz_{n+1}}\right|
-\frac{V_{\btheta^*_1,x^{n+1},T_{n+1},\bz_{n+1}}}{2}
\notag\\
&\leq U_{\btheta^*_1,x^{n+1},T_{n+1},\bz_{n+1}}-\frac{V_{\btheta^*_1,x^{n+1},T_{n+1},\bz_{n+1}}}{2}\notag\\
&\leq\log\hat f_{x^{n+1},T_{n+1},\bz_{n+1}}(X_{n+1})\notag\\
&\leq U_{\btheta^{**}_1,x^{n+1},T_{n+1},\bz_{n+1}}-\frac{V_{\btheta^{**}_1,x^{n+1},T_{n+1},\bz_{n+1}}}{2}\notag\\
&\leq \left|U_{\btheta^{**}_1,x^{n+1},T_{n+1},\bz_{n+1}}\right|+\frac{V_{\btheta^{**}_1,x^{n+1},T_{n+1},\bz_{n+1}}}{2}.\notag
\end{align}
Proceeding in the same way as in the $iid$ case, and exploiting Lemma \ref{lemma:moment_existence_non_iid}, we obtain
\begin{equation}
\underset{n}{\sup}~E\left(\log \hat f_{x^{n+1},T_{n+1},\bz_{n+1}}(X_{n+1})\right)^2<\infty.
\label{eq:var_num_finite_non_iid}
\end{equation}
Thus, as in the $iid$ set-up, (\ref{eq:convergent_series2}) follows from 
(\ref{eq:var_deno_finite_non_iid}) and (\ref{eq:var_num_finite_non_iid}). 

We formalize the above arguments in the form of a theorem in our non-$iid$ $SDE$ set-up. 
\begin{theorem}
\label{theorem:bf_consistency_non_iid}
Assume the non-$iid$ $SDE$ set-up and conditions 
(H1$^\prime$) -- (H6$^\prime$).
Then (\ref{eq:bf_consistency_non_iid}) holds.
\end{theorem}
As in the previous cases, the following corollary provides asymptotic comparison between two models
using Bayes factor in the non-$iid$ $SDE$ set-up.
\begin{corollary}
\label{corollary:sde_non_iid}
For $j=1,2$, let $R_{jn}(\btheta_j)=\prod_{i=1}^n\frac{f_{\btheta_j,x^i,T_i,\bz_i}(X_i)}{f_{\btheta_0,x^i,T_i,\bz_i}(X_i)}$, where
$\btheta_1$ and $\btheta_2$ are two different finite sets of parameters, perhaps with different dimensionalities, 
associated with the two models
to be compared. For $j=1,2$, let
\[
I_{jn}=\int R_{jn}(\btheta_j)\pi_j(d\btheta_j),
\]
where $\pi_j$ is the prior on $\btheta_j$. 
Let $B_n=I_{1n}/I_{2n}$ as before. 
Assume the non-$iid$ $SDE$ set-up and suppose that both
the models satisfy (H1$^\prime$) -- (H6$^\prime$), and
have the Kullback-Leibler property 
with $\delta=\delta_1$ and $\delta=\delta_2$ respectively. 
Then 
\[
n^{-1}\log B_n\rightarrow \delta_2-\delta_1,
\]
almost surely.
\end{corollary}

\section{Simulation studies}
\label{sec:simulated_data}

\subsection{Covariate selection when $n=15$, $T=1$}
\label{subsec:N15T1}
We demonstrate with simulation study the finite sample analogue of Bayes factor analysis
as $n\rightarrow\infty$ and $T$ is fixed. In this regard, we consider $n=15$ individuals, where the $i$-th one is modeled by
\begin{equation}
dX_i(t)=(\xi_1+\xi_2z_1(t)+\xi_3z_2(t)+\xi_4z_3(t))(\xi_5+\xi_6X_i(t))dt+\sigma_i dW_i(t),
\label{eq:sde2_appl}
\end{equation}
for $i=1,\cdots,15$. We fix our diffusion coefficients as $\sigma_{i+1}=\sigma_i +5$ for $i=1\cdots,14$ 
where $\sigma_1=10$. We consider the initial value $X(0)=0$ and the time interval $[0, T]$ with $T=1$. 

To achieve numerical stability of the marginal likelihood corresponding to each data 
we choose the true values of $\xi_i$; $i=1,\ldots,6$ as follows:
$\xi_{i}\stackrel{iid}{\sim} N(\mu_i,0.001^2)$, where 
$\mu_i\stackrel{iid}{\sim} N(0,1)$. This is not to be interpreted as the prior; this
is just a means to set the true values of the parameters of the data-generating model.
\\[2mm]
We assume that the time dependent covariates $z_i(t)$ satisfy the following $SDE$s 
\begin{align}
dz_1(t)=&(\tilde\theta_{1}+\tilde\theta_{2}z_1(t))dt+ dW_1(t)\notag\\
dz_2(t)=&\tilde\theta_{3}dt+ dW_2(t)\notag\\
dz_3(t)=&\tilde\theta_{4}z_3(t))dt+ dW_3(t),
\label{eq:covariate_appl}
\end{align}
where $W_i(\cdot)$; $i=1,2,3$, are independent Wiener processes, and 
$\tilde\theta_{i}\stackrel{iid}{\sim} N(0,0.01^2)$ for $i=1,\cdots,4$. 

We obtain the covariates by first simulating $\tilde\theta_{i}\stackrel{iid}{\sim} N(0,0.01^2)$ for $i=1,\cdots,4$,
fixing the values, and then by simulating the covariates using the $SDE$s (\ref{eq:covariate_appl}) 
by discretizing the time interval $[0,1]$ into $500$ equispaced time points. 
In all our applications we have standardized the covariates over time so that they have zero means
and unit variances.

Once the covariates are thus obtained, we assume that the data are generated from the (true) model where all the covariates
are present. For the true values of the parameters, we simulated $(\xi_1,\ldots,\xi_6)$ from
the prior and treated the obtained values as the true set of parameters $\btheta_0$.
We then generated the data using (\ref{eq:sde2_appl}) by discretizing the time interval $[0,1]$ 
into $500$ equispaced time points. 
%

As we have three covariates so we will have $2^3=8$ different models. Denoting
a model by the presence and absence of the respective covariates, it then is the case that 
$(1,1,1)$ is the true, data-generating model, while $(0,0,0)$, $(0,0,1)$, $(0,1,0)$, $(0,1,1)$, $(1,0,0)$,
$(1,0,1)$, and $(1,1,0)$ are the other $7$ possible models. 

\subsubsection{Case 1: the true parameter set $\btheta_0$ is fixed}
\label{subsubsec:case1_appl}
{\bf Prior on $\btheta$}
\\[2mm]
For the prior $\pi$ on $\btheta$, we first 
obtain the maximum likelihood estimator ($MLE$) of $\btheta$ using simulated annealing
(see, for example, \ctn{Liu01}, \ctn{Robert04}), and
consider a normal prior where the mean is the $MLE$ of $\xi_i$ for $i=1,\ldots,6$ 
and the variance is $0.8^2\mathbb I_6$, $\mathbb I_6$ being the $6$-dimensional identity matrix. 
As will be seen, this results in consistent model selection using Bayes factor.
\\[2mm]
{\bf Form of the Bayes factor}
\\[2mm]
In this case the related Bayes factor has the form
\begin{equation*}
I_n=\int \prod_{i=1}^n\frac {f_{i,\btheta_1}(X_i)}{f_{i,\btheta_0}(X_i)}\pi (d\btheta_1),
\label{eq:bf_T_appl}
\end{equation*}
where $\btheta_0=(\xi_{0,1},\xi_{0,2},\xi_{0,3},\xi_{0,4},\xi_{0,5},\xi_{0,6})$ is the true parameter set 
and $\btheta_1=(\xi_1,\xi_2,\xi_3,\xi_4,\xi_5,\xi_6)$ is the unknown set of parameters corresponding to 
any other model. 
Table \ref{table:values} describes the results of our Bayes factor analyses. 
\begin{table}[h]
\centering
\caption{Bayes factor results}
\label{table:values}
\begin{tabular}{|c||c|}
\hline
Model & $\frac{1}{15}\log I_{15}$\\
\hline
$(0,0,0)$ & -3.25214\\ 
$(0,0,1)$ & -1.39209\\
$(0,1,0)$ & -3.31954\\
$(0,1,1)$ & -1.11729\\
$(1,0,0)$ & -3.40378\\
$(1,0,1)$ & -1.22529\\
$(1,1,0)$ & -3.46790\\
\hline
\end{tabular}
\end{table}
It is clear from the 7 values of the table that the correct model $(1,1,1)$ is always preferred.

\subsubsection{Case 2: the parameter set $\btheta_0$ is random and has the prior distribution $\pi$}
\label{subsubsec:case2_appl}
We consider the same form of the prior $\pi$ as in 
Section \ref{subsubsec:case1_appl}, but with
variance $0.1^2\mathbb I_6$. The smaller variance compared to that in Case 1 attempts to somewhat 
compensate, in essence, for the lack of precise information about the true parameter values. 

In this case we calculate the marginal log-likelihood of the $8$ possible models as
$$\ell_i=\frac{1}{15}\log\int \prod_{i=1}^nf_{i,\btheta_1}(X_i)\pi (d\btheta_1);~i=1,\ldots,8,$$
with $\ell_8$ corresponding to the true model.
Table \ref{table:deltavalues} shows that $\ell_8$ 
is the highest. This clearly implies that 
the Bayes factor consistently selects the correct set of covariates even though the parameters of the true model are not fixed.
\begin{table}[h]
\centering
\caption{Values of $\frac{1}{15}\times$ marginal log-likelihoods}
\label{table:deltavalues}
\begin{tabular}{|c||c|}
\hline
Model & $\ell_i$\\
\hline
$(0,0,0)$ & 2.42430\\ 
$(0,0,1)$ & 4.29608\\
$(0,1,0)$ & 1.75213\\
$(0,1,1)$ & 4.84717\\
$(1,0,0)$ & 1.56242\\
$(1,0,1)$ & 4.92628\\
$(1,1,0)$ & 0.47111\\
$(1,1,1)$ & 5.84665 ($\mbox{true model}$)\\
\hline
\end{tabular}
\end{table}

\section{Summary and conclusion}
\label{sec:conclusion}
In this article we have investigated the asymptotic theory of Bayes factors when the models are associated with
systems of $SDE$'s consisting of sets of time-dependent covariates. The model selection problem we consider
encompasses appropriate selection of a subset of covariates, as well as appropriate selection of the part
of the drift function that does not involve covariates.
Such an undertaking, according to our knowledge, is a first-time effort which did not 
hitherto take place in the literature. 

We have established almost sure exponential convergence of the Bayes factor when the time domains remain bounded
but the number of individuals tend to infinity, in both $iid$ and non-$iid$ cases.  
In the non-$iid$ context, we
proposed and proved general results on Bayes factor asymptotics, which 
should be of independent interest.

Our simulation studies demonstrate that Bayes factor is a reliable criterion even in non-asymptotic
situations for capturing the correct set of covariates in our $SDE$ set-ups.

Note that our theory for non-$iid$ situations 
readily extends to model comparison 
problems when one of the models is associated
with an $iid$ system of $SDE$'s and another with a non-$iid$ system of $SDE$'s.
For instance, if the true model is associated with an $iid$ system, then $f_{0i}\equiv f_0\equiv f_{\btheta_0}$, and
the rest of the theory remains the same as our non-$iid$ theory of Bayes factors. The case when the other model is associated 
with an $iid$ system is analogous.

\section*{Acknowledgments}
The first author gratefully acknowledges her CSIR Fellowship, Govt. of India.

\newpage

\renewcommand\thefigure{S-\arabic{figure}}
\renewcommand\thetable{S-\arabic{table}}
\renewcommand\thesection{S-\arabic{section}}

\setcounter{section}{0}
\setcounter{figure}{0}
\setcounter{table}{0}

\begin{center}
{\bf \Large Supplementary Material}
\end{center}

Throughout, we refer to our main manuscript \ctn{Maitra15} as MB. 

\section{Proof of Lemma 1 of MB}
\label{sec:proof_lemma_1}
We first consider $k=1$. 
Note that due to assumption (H4$^\prime$), 
$$E_{\btheta_0}\left(V_{i,\btheta_1}\right)\leq 
T_i K_{\bbeta_1} \left(\underset{s\in [0,T_i]}{\sup}~\phi^2_{i,\bxi_j}(s)\right)
\left(1+\underset{s\in [0,T_i]}{\sup}~E_{\btheta_0}\left(X^2_i(s)\right)+\|\bbeta_1\|^2\right)<\infty,$$
since, by Proposition 1 of \ctn{Maud12}, $\underset{s\in [0,T_i]}{\sup}~E_{\btheta_0}\left(X^{2\ell}_i(s)\right)<\infty$,
for $\ell\geq 1$, and since $ \left(\underset{s\in [0,T_i]}{\sup}~\phi^2_{i,\bxi_j}(s)\right)$ is bounded
above due to continuity of $g_l$; $l=1,\ldots,p$. Hence, (3.6) of MB 
holds.

Now observe that due to Cauchy-Schwartz and (H4$^\prime$) 
\begin{align}
E_{\btheta_0}\left(V_{i,\btheta_0,\btheta_j}\right)
&=\int_0^{T_i}E_{\btheta_0}\left(\frac{\phi_{i,\bxi_0}(s)\phi_{i,\bxi_j}(s)
b_{\bbeta_j}(X_i(s))b_{\bbeta_0}(X_i(s))}{\sigma^2(X_i(s))}\right)ds\notag\\
&\leq \int_0^{T_i}\left[E_{\btheta_0}\left(\frac{\phi^2_{i,\bxi_j}(s)b^2_{\bbeta_j}(X_i(s))}
{\sigma^2(X_i(s))}\right)\right]^{\frac{1}{2}}
\times \left[E_{\btheta_0}\left(\frac{\phi^2_{i,\bxi_0}(s)b^2_{\bbeta_0}(X_i(s))}
{\sigma^2(X_i(s))}\right)\right]^{\frac{1}{2}}ds\notag\\
&\leq T_i\left(\underset{s\in [0,T_i]}{\sup}~\phi^2_{i,\bxi_j}(s)\right)^{\frac{1}{2}}K^{\frac{1}{2}}_{\bbeta_j} 
\left(1+\underset{s\in [0,T_i]}{\sup}~
E_{\btheta_0}\left(X^2_i(s)\right)+\|\bbeta_j\|^2\right)^{\frac{1}{2}}\notag\\
&\quad\quad\quad\times K^{\frac{1}{2}}_{\bbeta_0}\left(\underset{s\in [0,T_i]}{\sup}~\phi^2_{i,\bxi_0}(s)\right)^{\frac{1}{2}}
\left(1+\underset{s\in [0,T_i]}{\sup}~E_{\btheta_0}\left(X^2_i(s)\right)+\|\bbeta_0\|^2\right)^{\frac{1}{2}}
\notag\\
&<\infty,\notag
\end{align}
by Proposition 1 of \ctn{Maud12} and boundedness of 
$\left(\underset{s\in [0,T_i]}{\sup}~\phi^2_{i,\bxi_j}(s)\right)$. 
Hence, (3.7) of MB 
holds.
Also note that since
$E_{\btheta_0}\left(U_{i,\btheta_j}\right)=E_{\btheta_0}\left(V_{i,\btheta_0,\btheta_j}\right)$
by (2.11) and (2.13), 
(3.5) 
is implied by (3.7) of MB. 
To see that the moments are continuous in $\theta_1$, 
let $\left\{\btheta^{(m)}_1\right\}_{m=1}^{\infty}$ 
be a sequence converging to $\tilde\btheta_1$ as $m\rightarrow\infty$. 
Due to (H3$^\prime$), 
$$\frac{\phi^2_{i,\bxi^{(m)}_1}(s)b^2_{\bbeta^{(m)}_1}(X_i(s))}{\sigma^2(X_i(s))}\rightarrow 
\frac{\phi^2_{i,\tilde\bxi_1}(s)b^2_{\tilde\bbeta_1}(X_i(s))}{\sigma^2(X_i(s))}$$
and
$$\frac{\phi_{i,\bxi^{(m)}_1}(s)\phi_{i,\bxi_0}(s)b_{\bbeta^{(m)}_1}(X_i(s))b_{\bbeta_0}(X_i(s))}{\sigma^2(X_i(s))}\rightarrow 
\frac{\phi_{i,\tilde\bxi_1}(s)\phi_{i,\bxi_0}(s)b_{\tilde\bbeta_1}(X_i(s))b_{\bbeta_0}(X_i(s))}{\sigma^2(X_i(s))},$$
as $m\rightarrow\infty$, for any given sample path $\left\{X_i(s):s\in [0,T_i]\right\}$.
Assumption 
(H4$^\prime$) implies that 
$\frac{\phi^2_{i,\bxi^{(m)}_1}(s)b^2_{\bbeta^{(m)}_1}(X_i(s))}{\sigma^2(X_i(s))}$
is dominated by 
$\underset{\bxi_1\in\mathfrak \bGamma,s\in [0,T_i]}{\sup}~ \phi^2_{i,\bxi_1}(s)\times\underset{\bbeta_1\in\mathfrak B}{\sup}~ K_{\bbeta_1} 
\left(1+\underset{s\in [0,T_i]}{\sup}~\left[X_i(s)\right]^2+\underset{\bbeta_1\in\mathfrak B}{\sup}~\|\bbeta_1\|^2\right)$.
Since $X_i(s)$ is continuous on $[0,T_i]$, (guaranteed by (H2$^\prime$); see \ctn{Maud12}), it follows that
$\int_0^{T_i}\left[X_i(s)\right]^2ds<\infty$, which, in turn guarantees, in conjunction with
compactness of $\mathfrak B$ and $\bGamma$, that the upper bound is integrable. Hence, 
$V_{i,\btheta^{(m)}_1}\rightarrow V_{i,\tilde\btheta_1}$, almost surely.
Now, for all $m\geq 1$, 
\begin{equation*}
V_{i,\btheta^{(m)}_1}< T_i\left(\underset{\bxi_1\in\mathfrak \bGamma,s\in [0,T_i]}{\sup}~ \phi^2_{i,\bxi_1}(s)\right)
\times\left(\underset{\bbeta_1\in\mathfrak B}{\sup}~ K_{\bbeta_1}\right)
\times\left(1+\underset{s\in [0,T_i]}{\sup}~\left[X_i(s)\right]^2+\underset{\bbeta_1\in\mathfrak B}{\sup}~\|\bbeta_1\|^2\right).
\label{eq:v_bound}
\end{equation*}
Since $E_{\btheta_0}\left(\underset{s\in [0,T_i]}{\sup}~\left[X_i(s)\right]^2\right)<\infty$ 
by (3.1) of MB, 
it follows that
$E_{\btheta_0}\left(V_{i,\btheta^{(m)}_1}\right)
\rightarrow E_{\btheta_0}\left(V_{i,\tilde\btheta_1}\right)$,
as $\left(\btheta^{(m)}_1\right)\rightarrow \left(\tilde\btheta_1\right)$. 
Hence, $E_{\btheta_0}\left(V_{i,\btheta_1}\right)$
is continuous in $\btheta_1$. 

In the case of $V_{i,\btheta_0,\btheta_1}$, the relevant quantity
$\frac{\phi_{i,\bxi^{(m)}_1}\phi_{i,\bxi_0}b_{\bbeta^{(m)}_1}(X_i(s))b_{\bbeta_0}(X_i(s))}{\sigma^2(X_i(s))}$ 
is dominated by 
the continuous function (hence integrable on $[0,T_i]$) 
\begin{align}
&\left(\underset{\bxi_1\in\bGamma,s\in [0,T_i]}{\sup}~\left|\phi_{i,\bxi_1}(s)\right|\right)
\times\left(\underset{\bbeta_1\in\mathfrak B}{\sup}~ K^{\frac{1}{2}}_{\bbeta_1}\right)\times
\left(1+\underset{s\in [0,T_i]}{\sup}~\left[X_i(s)\right]^2+\underset{\bbeta_1\in\mathfrak B}{\sup}~\|\bbeta_1\|^2\right)\notag\\
&\quad\quad\quad\times
\left(\underset{s\in [0,T_i]}{\sup}~\left|\phi_{i,\bxi_0}(s)\right|\right)\times
K^{\frac{1}{2}}_{\bbeta_0}
\left(1+\underset{s\in [0,T_i]}{\sup}~\left[X_i(s)\right]^2+\|\bbeta_0\|^2\right),\notag
\end{align}
which ensures $V_{i,\btheta_0,\btheta^{(m)}_1}\rightarrow 
V_{i,\btheta_0,\tilde\btheta_1}$, almost surely.
Using the above bound for
\\
$\frac{\phi_{i,\bxi^{(m)}_1}\phi_{i,\bxi_0}b_{\bbeta^{(m)}_1}(X_i(s))b_{\bbeta_0}(X_i(s))}{\sigma^2(X_i(s))}$, 
it is seen that 
$$V_{i,\btheta_0,\btheta^{(m)}_1}< T_iK_1\left(\underset{s\in [0,T_i]}{\sup}~\left[X_i(s)\right]^4+K_2
\underset{s\in [0,T_i]}{\sup}~\left[X_i(s)\right]^2+K_3\right),$$
for appropriate positive constants $K_1,K_2,K_3$, so that (3.1) of MB 
for $k=4$, guarantees that
$E_{\btheta_0}\left(V_{i,\btheta_0,\btheta^{(m)}_1}\right)\rightarrow 
E_{\btheta_0}\left(V_{i,\btheta_0,\tilde\btheta_1}\right)$,
as $\btheta^{(m)}_1\rightarrow \tilde\btheta_1$. 
This shows that $E_{\btheta_0}\left(V_{i,\btheta_0,\btheta_1}\right)$
is continuous as well.
Since $E_{\btheta_0}\left(U_{i,\btheta_1}\right)
=E_{\btheta_0}\left(V_{i,\btheta_0,\btheta_1}\right)$, 
it follows that $E_{\btheta_0}\left(U_{i,\btheta_1}\right)$ is continuous in $\btheta_1$.

We now consider $k\geq 2$.
Note that, due to (H4$^\prime$), and the inequality $(a+b)^k\leq 2^{k-1}(|a|^{k}+|b|^{k})$ for $k\geq 2$ and any $a,b$, 
\begin{align}
E_{\btheta_0}\left(V_{i,\btheta_j}\right)^k
&\leq \left(\underset{s\in[0,T_i]}{\sup}~\left|\phi_{i,\bxi_j}(s)\right|^k\right)
2^{k-1}T_i^kK^k_{\bbeta_j} \left(1+\|\bbeta_j\|^2\right)^k\notag\\
&\quad\quad+\left(\underset{s\in[0,T_i]}{\sup}~\left|\phi_{i,\bxi_j(s)}\right|^k\right)2^{k-1}T_i^kK^k_{\bbeta_j}
E\left(\underset{s\in[0,T_i]}{\sup}~\left[X_i(s)\right]^{2k}\right).\notag
\end{align}
Since $E\left(\underset{s\in[0,T_i]}{\sup}~\left[X_i(s)\right]^{2k}\right)<\infty$ due to (3.1) of MB, 
and because $K_{\bbeta_j}$, $\|\bbeta_j\|$ are continuous in compact $\mathfrak B$, and 
$\left(\underset{s\in[0,T_i]}{\sup}~\left|\phi_{i,\bxi_j(s)}\right|^k\right)$ is continuous in compact
$\bGamma$, it holds that $E_{\btheta_0}\left(V_{i,\btheta_j}\right)^k<\infty$. In a similar manner it can be shown that
$E_{\btheta_0}\left(V_{i,\btheta_0,\btheta_1}\right)^k<\infty$. Thus, (3.7) of MB 
follows.

To see that (3.5) of MB 
holds, note that, due to (2.11) of MB 
and 
$(a+b)^k\leq 2^{k-1}(|a|^{k}+|b|^{k})$,
\begin{equation}
E_{\btheta_0}\left(U_{i,\btheta_j}\right)^k\leq 2^{k-1}E_{\btheta_0}\left(V_{i,\btheta_0,\btheta_j}\right)^k
+2^{k-1}E_{\btheta_0}\left(\int_0^{T_i}
\frac{\phi_{i,\bxi_j}(s)b_{\bbeta_j}\left(X_i(s)\right)}{\sigma
\left(X_i(s)\right)}dW_i(s)\right)^k. 
\label{eq:u_v_relation2}
\end{equation}
Since, due to (H4$^\prime$), (3.1) of MB 
and continuity of $\phi_{i,\bxi_j}$ on compact spaces, 
$$E_{\btheta_0}\left(\int_0^T\left|\frac{\phi_{i,\bxi_j}(s)b_{\bbeta_j}(X_i(s))}
{\sigma(X_i(s))}\right|^kds\right)<\infty,$$
Theorem 7.1 of \ctn{Mao11} (page 39) shows that 
\begin{equation}
E_{\btheta_0}
\left(\left|\int_0^{T_i}\frac{\phi_{i,\bxi_j}(s)b_{\bbeta_j}(X_i(s))}
{\sigma(X_i(s))}dW(s)\right|^k\right)
\leq\left(\frac{k(k-1)}{2}\right)^{\frac{k}{2}}T_i^{\frac{k-2}{2}}
E_{\btheta_0}\left(\int_0^{T_i}\left|\frac{\phi_{i,\bxi_j}(s)b_{\bbeta_j}(X_i(s))}
{\sigma(X_i(s))}\right|^kds\right).
\label{eq:moment2}
\end{equation}
Combining (\ref{eq:u_v_relation2}) with (\ref{eq:moment2}) and the result 
$E_{\btheta_0}\left(V_{i,\btheta_0,\btheta_1}\right)^k<\infty$, it follows that
$E_{\btheta_0}\left(U_{i,\btheta_j}\right)^k<\infty$.


As regards continuity of the moments for $k\geq 2$, first note that in the context of $k=1$, 
we have shown almost sure continuity of $V_{i,\btheta_1}$ with respect to $\btheta_1$.
Hence, $V^k_{i,\btheta_1}$ is almost surely continuous with respect to $\btheta_1$. That is, 
$\btheta^{(m)}_1\rightarrow\tilde\btheta_1$ 
implies $V^k_{i,\btheta^{(m)}_1}\rightarrow V^k_{i,\tilde\btheta_1}$,
almost surely. Once again, dominated convergence theorem allows us to conclude that
$E_{\btheta_0}\left(V_{i,\btheta^{(m)}_1}\right)^k\rightarrow 
E_{\btheta_0}\left(V_{i,\tilde\btheta_1}\right)^k$,
implying continuity of $E_{\btheta_0}\left(V_{i,\btheta_1}\right)^k$ with respect to $\btheta_1$.
Similarly, it is easy to see that $E_{\btheta_0}\left(V_{i,\btheta_0,\btheta_1}\right)^k$ 
is continuous with respect to $\btheta_1$. 
To see continuity of $E_{\btheta_0}\left(U_{i,\btheta_1}\right)^k$, first note that
$$E_{\btheta_0}\left[\int_0^{T_i}\left(\frac{\phi_{i,\bxi^{(m)}_1}(s)b_{\bbeta^{(m)}_1}(X_i(s))}
{\sigma(X_i(s))}
-\frac{\phi_{i,\tilde\bxi_1}(s)b_{\tilde\bbeta_1}(X_i(s))}
{\sigma(X_i(s))}\right)^2ds\right]\rightarrow 0,$$
as $m\rightarrow\infty$. The result follows as before by first noting pointwise convergence, and then
using (H4$^\prime$) and then (3.1) of MB, 
along with (H1$^\prime$) and boundedness
of $\phi_{i,\bxi^{(m)}_1}$. By It\^{o} isometry it holds that
$$E_{\btheta_0}\left[\int_0^{T_i}\frac{\phi_{i,\bxi^{(m)}_1}(s)b_{\bbeta^{(m)}_1}(X_i(s))}
{\sigma(X_i(s))}dW_i(s)
-\int_0^{T_i}\frac{\phi_{i,\tilde\bxi_1}(s)b_{\tilde\bbeta_1}(X_i(s))}{\sigma(X_i(s))}dW_i(s)\right]^2
\rightarrow 0.$$
Hence, 
$$\int_0^{T_i}\frac{\phi_{i,\bxi^{(m)}_1}(s)b_{\bbeta^{(m)}_1}(X_i(s))}{\sigma(X_i(s))}dW_i(s)
\rightarrow \int_0^{T_i}\frac{\phi_{i,\tilde\bxi_1}(s)b_{\tilde\bbeta_1}(X_i(s))}{\sigma(X_i(s))}dW_i(s)$$
in probability, as $m\rightarrow\infty$. Since $V_{i,\btheta_0,\btheta_1^{(m)}}
\rightarrow V_{i,\btheta_0,\tilde\btheta_1}$ 
almost surely as
$m\rightarrow\infty$, it follows from (2.11) of MB 
that 
$U_{i,\btheta^{(m)}_1}\rightarrow U_{i,\tilde\btheta_1}$ in probability, so that
$U^k_{i,\btheta^{(m)}_1}\rightarrow U^k_{i,\tilde\btheta_1}$ in probability. Using
(H4$^\prime$), (3.1) of MB 
and (H1$^\prime$), it is easily seen, using the same methods
associated with (\ref{eq:u_v_relation2}) and (\ref{eq:moment2}), that
$\underset{m}{\sup}~E_{\btheta_0}\left(U_{i,\btheta^{(m)}_1}\right)^{2k}<\infty$, proving that
$\left\{U^k_{i,\btheta^{(m)}_1}\right\}_{m=1}^{\infty}$ is uniformly integrable.
Hence, $E_{\btheta_0}\left(U_{i,\btheta^{(m)}_1}\right)^k\rightarrow 
E_{\btheta_0}\left(U_{i,\tilde\btheta_1}\right)^k$.
In other words, $E_{\btheta_0}\left(U_{i,\btheta_1}\right)^k$ is continuous in $\btheta_1$.

\section{Proof of Theorem 7 of MB}
\label{sec:proof_theorem_7}
Let us consider the martingale sequence
\[
S_N=\sum_{n=1}^{N}[\log(I_n/I_{n-1})+\mathcal K(f_{0n},\hat f_n)],
\]
which is a martingale because 
$E[\log(I_n/I_{n-1})|X_1,X_2,\ldots,X_{n-1}]=-\mathcal K(f_{0n},\hat f_n)$. 
Using the above it can be verified that if (5.8) of MB 
holds, implying
\[
\sum_{n=1}^{\infty} n^{-2} Var\left[\log\frac{I_n}{I_{n-1}}\right]<\infty,
\]
then $S_N/N\rightarrow 0$ almost surely. Therefore
\begin{equation}
N^{-1}\log I_N + N^{-1} \sum_{n=1}^{N}\mathcal K(f_{0n},\hat f_n) \rightarrow 0,
\label{eq:I_n_convergence1}
\end{equation}
almost surely, as $N\rightarrow\infty$.

Now consider $N^{-1}\sum_{i=1}^N\log\frac{f_{0i}(X_i)}{f_i(X_i)}$. If (5.7) of MB 
holds, implying
\[
\sum_{i=1}^{\infty} i^{-2} Var\left[\log\frac{f_{0i}(X_i)}{f_i(X_i)}\right]<\infty,
\]
then by Kolmogorov's strong law of large numbers in the independent but non-identical case,
\[
\frac{1}{N}\sum_{i=1}^N\log\frac{f_{0i}(X_i)}{f_i(X_i)}\rightarrow
\mathcal K^{\infty}(f_0,f)
\]
almost surely, as $N\rightarrow\infty$. 
Let  $\mathcal N_0(c)=\{f: \delta\leq\mathcal K^{\infty}(f_0,f)\leq\delta+c\}$, where $c>0$.
%
Now, note that,
%
\begin{align}
I_N &=\int\frac{\prod_{i=1}^N f_i(X_i)}{\prod_{i=1}^N f_{0i}(X_i)}\pi (\tilde {df})\notag\\
& \geq \int_{\mathcal N_0(c)} 
\exp\left( {\sum_{i=1}^N \log \frac {f_i(X_i)}{f_{0i}(X_i)}}\right)\pi (\tilde {df})\notag\\
& =\int_{\mathcal N_0(c)} 
\exp\left(- {\sum_{i=1}^N \log \frac{f_{0i}(X_i)}{f_i(X_i)}}\right)\pi (\tilde {df}).\notag
\end{align}
By Jensen's inequality,
\begin{align}
\frac{1}{N}\log I_N&\geq -\int_{\mathcal N_0(c)}
\frac{1}{N}\left({\sum_{i=1}^N \log \frac{f_{0i}(X_i)}{f_i(X_i)}}\right)\pi (\tilde {df})
\label{eq:jensen}
\end{align}
The integrand on the right hand side converges to $\mathcal K^{\infty}(f_0,f)$, pointwise for every $f$, given any
sequence $\{X_i\}_{i=1}^{\infty}$ associated with the complement of some null set.
Since, for all such sequences, uniform integrability of the integrand is guaranteed by (5.7) of MB, 
it follows that
the right hand side of (\ref{eq:jensen}) converges to $-\int_{\mathcal N_0(c)}\mathcal K^{\infty}(f_0,f)\pi (\tilde {df})$ almost surely.
Hence, almost surely,
\begin{align}
\underset{N}{\liminf}~N^{-1}\log I_N&\geq -\int_{\mathcal N_0(c)}\mathcal K^{\infty}(f_0,f)\pi (\tilde {df})\notag\\
&\geq -(\delta+c)\pi\left(\mathcal N_0(c)\right)\notag\\
&\geq -(\delta+c).\notag
\end{align}
Since $c>0$ is arbitrary, it follows that
%
\begin{equation}
\underset{N}{\lim\inf}~ N^{-1}\log I_N \geq -\delta, 
\label{eq:liminf_I_N}
\end{equation}
almost surely.
%
Now, due to (5.2) of MB 
it follows that
$\mathcal K(f_{0n},f_n)\geq \delta$ for all $n$ with probability 1,
so that $\mathcal K(f_{0n},\hat f_n)=\mathcal K(f_{0n},\hat f_{n,A_n(\delta)})$,
where $A_n(\delta)$ is given by (5.5) of MB. 
By the $Q^*$ property 
it implies that  
\begin{equation*}
\underset{N}{\lim\inf}~ N^{-1} \sum_{n=1}^N \mathcal K(f_{0n},\hat f_n)\geq \delta.
\label{eq:liminf_K_average}
\end{equation*}
Hence, it follows from (\ref{eq:I_n_convergence1}) that
\begin{equation}
\underset{N}{\lim\sup}~ N^{-1} \log I_N \leq -\delta.
\label{eq:limsup_I_N}
\end{equation}
Combining (\ref{eq:liminf_I_N}) and (\ref{eq:limsup_I_N}) it follows that 
$$\underset{N\rightarrow\infty}{\lim}~N^{-1}\log I_N =-\delta,$$
almost surely.

\section{Proof of Lemma 10 of MB}
\label{sec:proof_lemma_9}

The proofs of (6.3) -- (6.5) of MB 
follow in the same way as the proofs of (3.5) -- (3.7) of MB,  
using compactness of $\mathfrak X$, $\mathfrak T$ 
and $\boldcal Z$ in addition to that of $\mathfrak B$ and $\bGamma$.

For the proofs of continuity of the moments, note that as in the $iid$ case,
uniform integrability is ensured by (H4$^\prime$), (3.1) of MB 
and compactness
of the sets $\mathfrak B$, $\bGamma$, $\mathfrak X$, $\mathfrak T$ and $\boldcal Z$. 
The rest of the proof is almost the same as the proof
of Theorem 5 of \ctn{Maitra14a}.

\section{Proof of Lemma 11 of MB}
\label{sec:proof_lemma_10}
For notational simplicity, let
\begin{align}
\breve V_{x^i,\bbeta_0}(s)&=\frac{b^2_{\bbeta_0}(X_i(s,x^i))}{\sigma^2(X_i(s,x^i))};\notag\\
\breve V_{x^i,\bbeta_0,\bbeta_j}(s)&=\frac{b_{\bbeta_j}(X_i(s,x^i))
b_{\bbeta_0}(X_i(s,x^i))}{\sigma^2(X_i(s,x^i))};\notag\\
\breve V_{x^i,\bbeta_1}(s)&=\frac{b^2_{\bbeta_1}(X_i(s,x^i))}{\sigma^2(X_i(s,x^i))}.\notag
\end{align}

Continuity of $\mathcal K(f_{\btheta_0,x,T,\bz},f_{\btheta_1,x,T,\bz})$ 
with respect to $x$ and $T$, the fact that $x^k\rightarrow x^{\infty}$ and 
$T_k\rightarrow T^{\infty}$ as $k\rightarrow\infty$, assumption (H5$^\prime$), and the dominated convergence theorem
together ensure that
\begin{align}
&\underset{n\rightarrow\infty}{\lim}~\frac{\sum_{k=1}^n\mathcal K(f_{\btheta_0,x^k,T_k,\bz_k},
f_{\btheta_1,x^k,T_k,\bz_k})}{n}\notag\\
&=\underset{n\rightarrow\infty}{\lim}~\frac{1}{n}
\sum_{k=1}^n\int_0^{T_k}\left\{\frac{{\phi_{k,\bxi_0}(s)}^2}{2}E_{\btheta_0}(\breve V_{x^k,\bbeta_0}(s))\right.
\notag\\
&\qquad\qquad\left.-\phi_{k,\bxi_0}(s)\phi_{k,\bxi_1}(s)E_{\btheta_0}(\breve V_{x^k,\bbeta_0,\bbeta_1}(s))
+\frac{{\phi_{k,\bxi_1}(s)}^2}{2}E_{\btheta_0}(\breve V_{x^k,\bbeta_1}(s))\right\}ds\notag\\
&=\underset{n\rightarrow\infty}{\lim}~\frac{1}{n}
\int_0^{T_k}\sum_{k=1}^n\left\{\frac{(\xi_{00}+\xi_{10}g_1(z_{k1}(s))
+\cdots+\xi_{p0}g_p(z_{kp}(s)))^2}{2}E_{\btheta_0}(\breve V_{x^k,\bbeta_0}(s))\right.\notag\\
&\quad\quad\left. -(\xi_{00}+\xi_{10}g_1(z_{k1}(s))+\cdots+\xi_{p0}g_p(z_{kp}(s)))(\xi_{01}
+\xi_{11}g_1(z_{k1}(s))+\cdots+\xi_{p1}g_p(z_{kp}(s)))E_{\btheta_0}(\breve V_{x^k,\bbeta_0,\bbeta_1}(s))\right.\notag\\
&\quad\quad\left. +\frac{(\xi_{01}+\xi_{11}g_1(z_{k1}(s))+\cdots+
\xi_{p1}g_p(z_{kp}(s)))^2}{2}E_{\btheta_0}(\breve V_{x^k,\bbeta_1}(s))\right\}ds\notag\\
&=\int_0^{T^{\infty}}\left\{\left(\frac{\xi_{00}^2}{2}+\xi_{00}\sum_{l=1}^p\xi_{l0}c_l(s)
+\frac{1}{2}\sum_{l=1}^p\sum_{m=1}^p\xi_{l0}\xi_{m0}c_l(s)c_m(s)\right)
E_{\btheta_0}(\breve V_{x^{\infty},\bbeta_0}(s))\right.\notag\\
&\quad\quad\left.-\left(\xi_{00}\xi_{01}+\xi_{00}\sum_{l=1}^p\xi_{l1}c_l(s)+\xi_{01}\sum_{l=1}^p\xi_{l0}c_l(s)+
\sum_{l=1}^p\sum_{m=1}^p\xi_{l0}\xi_{m1}c_l(s)c_m(s)\right)E_{\btheta_0}(\breve V_{x^{\infty},\bbeta_0,\bbeta_1}(s))
\right.\notag\\
&\quad\quad\left.+\left(\frac{\xi_{01}^2}{2}+\xi_{01}\sum_{l=1}^p\xi_{l1}c_l(s)
+\frac{1}{2}\sum_{l=1}^p\sum_{m=1}^p\xi_{l1}\xi_{m1}c_l(s)c_m(s)\right)
E_{\btheta_0}(\breve V_{x^{\infty},\bbeta_1}(s))\right\}ds,
\label{eq:kl_limit_1}
\end{align}
which is the Kullback-Leibler distance between the models of the same form as
$f_{\btheta_0,x,T,\bz}$ and $f_{\btheta_1,x,T,\bz}$, but with $x$, $T$ and $g_l(z_l(s))$ replaced with
$x^{\infty}$, $T^{\infty}$ and $c_l(s)$.

\normalsize
\bibliographystyle{natbib}
\bibliography{irmcmc}

\end{document}